\documentclass[11pt, oneside]{article} 
\usepackage[margin=2cm]{geometry}                		
\usepackage[parfill]{parskip}    		
\usepackage{graphicx}			
\usepackage[utf8]{inputenc}
\usepackage[english]{babel} % English language/hyphenation
\usepackage{amsmath,amssymb,amsthm,bm} % Math packages
\usepackage[margin=2cm]{geometry}
\usepackage[color=yellow]{todonotes}
\usepackage{url}	
\usepackage{esint}
\usepackage[colorlinks]{hyperref}
\usepackage{enumerate}
\usepackage{outlines}[enumerate]
\usepackage{framed,comment}
\usepackage{upgreek}
\usepackage{pgfplots}
\usepackage{float}
\usepackage{tikz,tikz-cd}
\usepackage{rotating}
\usepackage{graphicx}
\usepackage{caption}
\usepackage{multicol}
\usepackage{cancel}
\usepackage{subcaption}
\usepackage[title]{appendix}
\usepackage{tikz-cd}
\usepackage{pgf, pgffor}
\extrafloats{100}
\numberwithin{equation}{section}

\newtheorem{remark}{Remark}[section]

\newcommand{\mb}[1]{\mathbf{#1}}
\newcommand{\mbb}[1]{\mathbb{#1}}

\newcommand{\mc}[1]{\mathcal{#1}}
\newcommand{\bs}[1]{\boldsymbol{#1}}

\newcommand{\scp}[2]{\left<#1\,,\,#2\right>}
\newcommand{\bracket}[2]{\left<#1\,|\,#2\right>}
\newcommand{\ad}{\operatorname{ad}}

\def\p{{\partial}}

\def\rmd{{\color{red}{\rm d}}}
\def\bk{{\mathbf{k}}}

\def\bM{{\mathbf{M}}}
\def\bR{{\mathbf{R}}}
\def\bu{{\mathbf{u}}}
\def\bU{{\mathbf{U}}}
\def\bv{{\mathbf{v}}}
\def\bx{{\mathbf{x}}}

\def\Diff{\text{Diff}}

\def\p{\partial}

\begin{document}

\title{\textbf{Lagrangian reduction and wave mean flow interaction}}
\author{Darryl D. Holm, Ruiao Hu\footnote{Corresponding author, email: ruiao.hu15@imperial.ac.uk}, and Oliver D. Street\\
d.holm@imperial.ac.uk, ruiao.hu15@imperial.ac.uk, o.street18@imperial.ac.uk\\
Department of Mathematics, Imperial College London \\ SW7 2AZ, London, UK}
% 	\date{Key words: nonlinear water waves, free surface fluid dynamics, geometric mechanics}
\date{\bigskip
{\it \dots the most difficult task is
to think of workable examples that will reveal
something new \cite{Flaschka2015}.}}

\maketitle

\begin{abstract}
%takes a geometric viewpoint in addressing a central question arising in many fields of multiscale, %multiphysics, multicomponent hydrodynamics. The 
%This paper begins by asking the following question: 
How does one derive models of dynamical feedback effects in multiscale, multiphysics systems such as wave mean flow interaction (WMFI)? We shall address this question for hybrid dynamical systems, whose motion can be expressed as the composition of two or more Lie-group actions. Hybrid systems abound in fluid dynamics. Examples include: the dynamics of complex fluids such as liquid crystals; wind-driven waves propagating with the currents moving on the sea surface; turbulence modelling in fluids and plasmas; and classical-quantum hydrodynamic models in molecular chemistry. From among these examples, the motivating question in this paper is: How do wind-driven waves produce ocean surface currents? \\The paper first summarises the geometric mechanics approach for deriving hybrid models of multiscale, multiphysics motions in ideal fluid dynamics. It then illustrates this approach for WMFI in the examples of 3D WKB waves and 2D wave amplitudes governed by the nonlinear Schr\"odinger (NLS) equation propagating in the frame of motion of an ideal incompressible inhomogeneous Euler fluid flow. The results for these examples tell us that the mean flow in WMFI does not create waves. However, feedback in the opposite direction is possible, since 3D WKB and 2D NLS wave dynamics can indeed create circulatory mean flow.

%Although it is not discussed in detail, a stochastic geometric mechanics approach to this problem is also formulated, in preparation for further studies.
\end{abstract}

\tableofcontents

\section{Introduction}

{\bf Interaction of wind waves and ocean currents.} 
In the Iliad, one of Homer's verses describing air-sea interaction seems to hint that 
wind-driven waves convey an impulse of momentum into the sea \cite{IJtexts} 
\begin{quote}
like blasts of storming winds
striking the earth under Father Zeus’s thunder,
then with a roar slicing into the sea, whipping up
a crowd of surging waves across a booming ocean,
with lines of arching foam, one following another
\end{quote}

Modern geophysical fluid dynamics (GFD) would not disagree with Homer's simile for air-sea interaction. In particular, the well-known Craik-Leibovich (CL) theory of the role of Stokes drift in the creation of Langmuir circulations \cite{CL1976} and the Andrews-McIntyre theory of generalised Lagrangian mean (GLM) dynamics \cite{AM1978} each introduce a shift in the definition of total momentum due to fluctuating subsystem  dynamics and a corresponding non-inertial force on the mean fluid motion. 

In this paper, we use standard methods of geometric mechanics to formulate models of wave mean flow interaction (WMFI) of fluctuations on the Earth's mean sea surface that is based on boosting the fluctuation dynamics into the frame of the mean flow. We hope that such a model may become useful, for example, in the interpretation of satellite image data from the Surface Water and Ocean Topography (SWOT) satellite mission, which is the first satellite with the specific aim to measure fluctuations on the Earth's sea surface \cite{SWOT2023}. 

Our objective here is to construct WMFI dynamics as a \emph{hybrid} fluid theory based on symmetry reduction in an Euler-Poincar\'e variational principle for the nonlinear dynamics of a system of two fluid degrees of freedom \cite{HMR1998}. The mathematical theory formulated here is illustrated in a hybrid fluid theory reminiscent of Landau's two-fluid model of superfluid $He$-II as discussed, e.g., in \cite{Lin1963}. Just as with superfluids, the formulation of the theory in this paper involves transforming between the frames of motion of the two fluidic degrees of freedom. The role of the superfluid component of Landau's two-fluid $He$-II model in the WMFI model proposed here is played by the slowly varying complex amplitude of WKB wave equations, e.g., of the nonlinear Schr\"odinger (NLS) equation.  

% %%%%%%%%%%%%%%%%%%%%%%%%%%%%%%%%%%%%%%%%%%%%%%%%%%%%%%%%%%%%%%%%%%%%%%

In the absence of additional assumptions, the inverse problem of determining a three-dimensional fluid flow under gravity solely from observations of its two-dimensional surface flow and its propagating wave elevation field has no unique solution. Without attempting to discover the three-dimensional flow beneath the surface, though, one may still derive a mathematical model of some of the phenomena on the free surface via the implications of the kinematic boundary condition. Specifically, the kinematic boundary condition implies a composition of  horizontal flow and vertically oscillating wave elevation dynamics of the Lagrangian material parcels remaining on the surface. In this paper, we formulate the initial value problem for wave dynamics on the free surface of a three-dimensional fluid flow. This is done entirely in terms of surface phenomena, as the semi-direct composition of a two-dimensional area-preserving horizontal fluid flow map acting on the vertical wave elevation dynamics. The surface wave dynamics formulated here is derived via Hamilton's variational principle by using a Lagrangian comprising the difference of the fluid kinetic and potential energies, constrained by the kinematic boundary condition that the flow of material parcels remains on the surface of the fluid.

\subsection{Examples of hybrid models}

Hybrid systems often involve sequences of relative motions in which one degree of freedom evolves in the frame of motion of the previous one. Lewis Fry Richardson's ``whorls within whorls'' verse about the turbulence cascade describes the familiar situation in which big whorls, little whorls and lesser whorls interact  sequentially, one within the frame of motion of the one before, each feeling an additional reaction force from the change of frame. Plasma dynamics exemplifies another type of hybrid system, one in which Lagrangian charged particles interact with Eulerian electromagnetic fields. In this case, the Lorentz force on the charged fluid particles arises in the plasma fluid motion equation when the electromagnetic fields are Lorentz-transformed into the frame of the moving fluid. This type of reaction force due to a frame change can usually be attributed to a momentum shift associated with the change of frame for the electromagnetic field dynamics in the moving medium. 

{\bf Complex fluids.}
In a sequence of papers \cite{GBR09,GBRT2012,GBT2010,H2002} the geometric mechanics of perfect complex fluids was developed, culminating in a full description of the geometry underlying the classic Ericksen-Leslie and and Eringen theories of complex fluids\cite{GBRT2012}. The hybrid model approach we shall discuss in the present paper is consistent with these previous approaches. 

The next three hybrid models have the additional feature that the hybrid components of the degrees of freedom live in nested sets of physical spaces or phase spaces. 

{\bf Multiscale fluid turbulence models.}
The geometric hybrid approach also applies in the kinetic sweeping of microstructure in turbulence models \cite{HT2012}. The basic idea in these turbulence models is that the coarse-grained space contains the fine-grained space as a subgrid-scale degree of freedom. The fine-grained fluid dynamics are transported along the Lagrangian paths of the coarse-grained fluid dynamics by a composition of maps. Spatial averages over the evolution of the fine-grained fluid dynamics act back on the motion in the coarse-grained space and modify it. The back-reaction is calculated via the coarse-grained divergence of the Reynolds stress tensor for the coarse-grained fluid dynamics. The latter is defined by spatial averaging over the terms in the coarse-grained dynamics that feed back from the fluid dynamics in the fine-grained space, which is again parameterised by the coarse-grained coordinates by the composition of smooth invertible maps. 

{\bf Hybrid models of electromagnetic wave / fluid plasma interaction.} A natural candidate for hybrid models would be the electromagnetic wave / fluid plasma interaction. Examples of hybrid models of the geometric type considered here in plasma physics include: (i) ponderomotive coupling of microwaves and plasma in magnetic controlled fusion \cite{SKH1984,SKH1986}; Electro- and magneto- fluids \cite{H1986}; (iii) relativistic fluid plasma dynamics \cite{H1987}; and (iv) Vlasov--fluid hybrid plasma models, Holm and Kaufman \cite{HK1984}, Holm and Tronci 2010 \cite{HT2010}.

{\bf Classical--quantum mechanics.}
The coupling between classical and quantum degrees of freedoms has raised an outstanding question ever since the rise of quantum mechanics as a physical theory. How does one separate classical and quantum? How do they influence one another? Is there a back reaction? For example, is there something like Newton's Law of action and reaction when a classical measurement of a quantum property occurs?  A general model of classical--quantum back-reaction must be able to give consistent answers to the various quantum paradoxes. 

For example, the exact factorisation (EF) model  of quantum molecular chemistry is discussed from the viewpoint of the geometric mechanics approach in \cite{FHT2019,GBT2022,HRT2021,RT2021}. The EF model shares some similarities with the multiscale turbulence models in that two spatial scales are involved: one spatial scale for the slowly moving classical dynamics of the ions; and another spatial scale for the rapid quantum motion. The term ``exact factorisation" indicates that the total wave function is factorised into a classical wave function for the ions depending on one set of coordinates and a quantum wave function depending on a second set of coordinates whose motion relative to the first set of coordinates is determined by a composition of maps.

{\bf Image registration by LDM using the metamorphosis approach.}
Large deformation diffeomorphic matching methods (LDM) for image registration are based on optimal control theory, i.e., minimizing the sum of a kinetic energy metric, plus a penalty term. The former ensures that the template deformation by the diffeomorphism follows an optimal path, while the latter ensures an acceptable tolerance in image mismatch.
The \emph{metamorphosis approach} is a variant of LDM that parallels the present considerations, in allowing  dynamical templates, so that the evolution of the image template deviates from pure deformation \cite{TY2005,HTY2009}.

{\bf Wave mean flow interaction.} The hybrid description of WMFI in terms of two fluid fields is already standard in geophysical fluid dynamics (GFD) models. For example, the Craik-Leibovich (CL) approach \cite{CL1976} and the Generalised Lagrangian Mean (GLM) approach \cite{AM1978, GH1996} both introduce two types of fluid velocities, one for the mean flow and another for the fluctuations. See \cite{SFK2016} for a recent summary of the state of the art in Craik-Leibovich models.

{\bf The present work.}
In all of the hybrid models mentioned so far, a simple and universal property of transformation theory called the cotangent-lift momentum map plays a key role in describing the interactions among the various degrees of freedom in the hybrid dynamical system. The same property plays a key role in the theory developed here for the interaction of free-surface waves and the fluid currents which transport them. 

Thus, the present work extends the ongoing series of applications of geometric mechanics in multiscale, multiphysics continuum dynamics to the case of the interaction of fluid waves and currents. As mentioned earlier, we hope that restricting this approach to two spatial dimensions will contribute a useful method for data calibration and analysis of satellite observations of the ocean surface in the SWOT mission. In preparation for the data calibration, analysis, and assimilation aspects of the potential applications of this approach, we also include Appendix \ref{appendix: stochastic WCI} which formulates the stochastic versions of the deterministic WMFI equations treated in the main body of the paper that could be useful as a basis for SWOT data analysis. 

\paragraph{Plan of the paper.}
Section \ref{sec:2-Lagrangian reduction} shows the Lie group reduced variational path via Hamilton's principle for deriving hybrid fluid models.
In Section \ref{sec:3-NLS}, we introduce and discuss two examples of hybrid models. These hybrid fluid models are Eulerian wave elevation field equations governing the coupling of an Euler fluid to: (i) harmonic scalar wave field elevation oscillations; and (ii) complex scalar elevation field dynamics governed by the nonlinear Schr\"odinger (NLS) equation. The latter are called \emph{hybrid Euler-NLS equations}.
Section \ref{sec:4-Numerics} shows simulations of the hybrid Euler-NLS equations that verify the predictions of momentum exchange derived in the previous section.
Section \ref{sec:5-Conclusion} provides concluding remarks, as well as an outlook toward future work. Appendix \ref{appendix: stochastic WCI} proposes stochastic modifications of the present deterministic variantional theory and Appendix \ref{appendix:SHO} discusses an instructive elementary example in which the waves comprise a field of vertical simple harmonic oscillators.

\section{Lagrangian reduction}\label{sec:2-Lagrangian reduction}
We are dealing with physical problems that involve a subset of variables evolving in the frame of reference moving with an underlying fluid dynamical system. An example was given earlier of waves propagating in the frame of reference given by ocean currents \cite{Proceedings2022}. 
% Further examples include complex fluids \cite{}, the elastic pendulum \cite{}, or the double pendulum \cite{}.
In general, the dynamics of an order parameter breaks a symmetry that the system would have had in the abssence of said parameter. This problem may be described geometrically in the following way. Motivated by wave mean flow interactions (WMFI), within this section we will perform the calculations for the case of continuum dynamics, where the Lie group acting on the order parameters is taken to be the group of diffeomorphisms. We will therefore choose Lagrangians, group actions, and representations that are \emph{right} invariant. It should be noted that the theory presented in this section is general enough to apply for other dynamical systems whose behaviour can be described by the action of a Lie group on a configuration space. The abstract equations derived in this section may be derived similarly under the assumption that we are working with left actions and left invariant Lagrangians, whereupon the equations would have a slightly modified form in the same way as is known for standard Euler-Poincar\'e systems \cite{HMR1998}. 

The configuration space of fluid motion within a spatial domain\footnote{For our examples of WMFI dynamics, we will take dimension $n=3$ and $n=2$ for the examples in section \ref{sec:3-NLS}}, $\mathcal{D} \in \mathbb{R}^n$, is given by the diffeomorphism group, $G={\rm Diff}(\mathcal{D})$. That is, each element, $g\in G$, is a map from $\mathcal{D}$ to itself which takes a fluid particle at a position, $X\in\mathcal{D}$, at initial time $t=0$, to a position, $x = g_t(X)$, at the current time, $t$, with $g_0(X)=X$, so that $g_0=Id$. The time-parametrised curve of diffeomorphisms, $g_t\in G$, therefore governs the history of each fluid particle path within the domain. Thus, the fluid motion is described by the transitive action of $G$ on $\mathcal{D}$. In what follows, we will denote the corresponding Lie algebra by $\mathfrak{g}$, which for fluid motion is the space of vector fields, i.e.,  $\mathfrak{g}=\mathfrak{X}(\mathcal{D})$. 
 
For a $G$-invariant Lagrangian defined on the tangent bundle, $TG$, the equations of motion are given by the standard Euler-Poincar\'e theorem, which can be expressed on $G$, or in their reduced form on the dual of the Lie algebra, $\mathfrak{g}^*=\Lambda\otimes {\rm Den}(\mathcal{D})$, the 1-form densities on domain $\mathcal{D}$ in the case of fluids with $L^2$ pairing. The symmetry of this description can be broken by the presence of a \emph{parameter}, $a_0 \in V^*$, in a vector space $V^*$ which admits a right representation of $G$ acting on $V$. Note that the convention of choosing this space of parameters to be the dual space, $V^*$, is due to the development of semidirect product reduction theory for Hamiltonian mechanics \cite{MRW1984a,MRW1984b}, which occurred prior to its development on the Lagrangian side \cite{HMR1998}. For a Lagrangian, $L:TG\times V^* \to \mathbb{R}$, we may define a collection of Lagrangians, $L_{a_0}:TG \to \mathbb{R}$, parameterised smoothly by $a_0 \in V^*$. The induced right representation of $G$ on $V^*$ is given by push-forward, and thus $G$ acts on an element $a_0 \in V^*$ as
\begin{align}
a_t(x) = a_0(X)\cdot g_t^{-1} =: g_{t\,*}a_0(X)\,.
\label{advec-reln}
\end{align}
It is now apparent that the Lagrangian, $L_{a_0}$, is invariant under the action of the isotropy group\footnote{This group is sometimes referred to as the stabiliser of $a_0$.} of $a_0$, which is defined by
\begin{equation}
	G_{a_0} := \{ g \in G \ | \ g_{t\,*}a_0 = a_0 \} \,. 
\end{equation}
The Lagrangian is \emph{not} invariant under the action of $G$, and we thus say that the symmetry is \emph{broken}. The remaining invariance under the action of the isotropy group is the particle relabelling symmetry. Consider the group of cosets, $G \setminus G_{a_0}$. One may think equivalently of $G$ acting on $G \setminus G_{a_0}$
 or $V^*$ (where this equivalence follows from Theorem 21.18 in \cite{L2012}). Indeed, there is a correspondence between the parameter $a_0\in V^*$ and the coset $[e] = e\cdot G_{a_0} = G_{a_0}$, where $e$ is the identity element of $G$. If we instead choose a parameter $a_0'$, then since the action of $G$ on $V^*$ is transitive there exists a group element, $h\in G$, such that $a_0' = h_*a_0$. The isotropy group for $a_0'$ is the conjugate subgroup, $G_{a_0'} = h^{-1}G_{a_0}h$, and the coset spaces $G \setminus G_{a_0}$ and $G \setminus G_{a_0'}$ are hence isomorphic. Whilst the parameters, $a_0 \in V^*$, break the symmetry of the system, they may be chosen arbitrarily.
 
We note that $a_t$ as given in equation \eqref{advec-reln} is the solution of the advection equation, denoted as
\begin{align}
\p_ta+{\cal L}_u a=0\,,
\quad\hbox{with}\quad
u := \dot{g}_tg_t^{-1}
\,,
\label{advec-eqn}
\end{align}
where ${\cal L}_u$ denotes the Lie derivative with respect to the Eulerian velocity vector field, $u:=\dot{g}_tg_t^{-1}$.
The advection equation follows from the Lie chain rule for the push-forward $g_{t\,*}$ of the initial condition $a_0(X)$ by the time-dependent smooth invertible map $g_t$. Namely,
\begin{align}
\p_ta_t(x) = \p_t \big(g_{t\,*}a_0(X)\big) 
= -\, g_{t\,*} \big({\cal L}_{\dot{g}_tg_t^{-1}}a_0(X)\big)
= -\,{\cal L}_u a_t(x)
= -\,{\cal L}_u a(x,t)
\,.
\label{LieCR}
\end{align}

Imposing the advection relation in \eqref{advec-reln} in Hamilton's principle when the Lagrangian is right invariant under the action of $g_t$ yields the standard Euler-Poincar\'e theory for semidirect product Lie algebras, \cite{HMR1998}.

\paragraph{Additional dynamics.} Suppose further that we have an additional configuration space, $Q$, which represents (order) parameters with their own dynamics, and that we have a representation of the (free, transitive) group action of $G$ on $Q$.
Within this space we will find dynamics (e.g. waves) occurring within the frame of reference corresponding to the (fluid) motion on $\mathfrak{g}^*$.
The distinction between parameters in $V^*$ and $TQ$ becomes apparent in the variational formulation. Indeed, let us consider the general case in which the Lagrangian $L$ takes the form
\begin{align}
    L:T\left(G\times Q\right)\times V^* \rightarrow \mathbb{R}\,,
\end{align}
where $G$, $Q$, and $V^*$ are as defined above. We assume that $G$ acts on $T\left(G\times Q\right)\times V^*$ in the natural way on the right. We denote this right action using concatenation and tangent fibre notation $u_g$ at footpoint $g$ on the manifold $G$ as
\begin{align}
    (g, u_g, q, u_q, a)h = (gh, u_gh, qh, u_qh, ah)\,.
\end{align}
Invariance of the Lagrangian $L$ in Hamilton's principle under the right action of $G$ is written as
\begin{align}
    L(g, \dot{g}, q, \dot{q}, a_0) = L(gh, \dot{g}h, qh, \dot{q}h, a_0h) \,,
\end{align}
for all $h \in G$. Choosing $h = g^{-1}$, one defines the reduced Lagrangian as
\begin{align}
    L(e, \dot{g}g^{-1}, q g^{-1}, \dot{q}g^{-1}, a_0g^{-1}) =: \ell(u, n, \nu, a)
    \,,
\end{align}
with further notation $u:=\dot{g}g^{-1}$, $n = q g^{-1}$ and $\nu = \dot{q}g^{-1}$. The reduced Lagrangian $\ell$ is then associated to the quotient map,
\begin{align}
    T(G\times Q)\times V^* \rightarrow \mathfrak{g}\times TQ\times V^*\,.
\end{align}
We have thus formulated the reduced Euler--Poincar\'e variational principle,
\begin{align}
    0 = \delta S = \delta \int_{t_0}^{t_1} \ell(u, n, \nu, a)\,dt\,, \label{eq:metamorphsis var principle}
% \label{var-princ}
\end{align}
defined subject to the following constrained variations of $u, n, \nu$ and $a$, derived from their definitions,
\begin{align}
\begin{split}
    \delta u &= \p_t\eta - \ad_u \eta \,,\\
    \delta n &= w - \mathcal{L}_\eta n \,,\\
    \delta \nu &= \p_t w + \mathcal{L}_u w - \mathcal{L}_\eta \nu\,,\\
    \delta a &= -\mathcal{L}_\eta a\,,
\end{split} \label{eq:metamorphsis constrained variations}
\end{align}
where ${\rm ad}_u\eta = - [u,\eta]$, $\eta = \delta g g^{-1}$ and $w = \delta q g^{-1}$ are arbitrary and vanish at the endpoints in time, $t = t_0$ and $t=t_1$. Here, the Lie derivative
% \todo[inline]{OS: This isn't a Lie derivative yet, since $\eta$ is not a vector field in general. I guess we just call it $\eta q$? Not sure what notation to use instead of transpose of Lie derivative!\\
% DH: $\eta = \delta g g^{-1}$ is indeed a vector field for $g\in {\rm Diff}(Q)$. As you intended, right Ollie?} 
w.r.t to the vector field $\eta$ is denoted as $\mathcal{L}_\eta$. Note that the form of the arbitrary vector fields $\eta$ and $w$ is a consequence of the fact that our group action is from the right, which also dictates the form of $n$ and $\nu$. The dual Lie derivative operator, $\diamond$, is defined via nondegenerate pairings $\scp{\cdot}{\cdot}$ over $\mathfrak{g}$ and $T^*Q$ as
\begin{align}
    \scp{p}{\mathcal{L}_\eta q}_{T_q^*Q \times T_qQ} = \scp{-p\diamond q}{\eta}_{\mathfrak{g}}\,,
\label{def-diamond}
\end{align}
for all $(p,q)\in T^*Q$ and $\eta\in\mathfrak{g}$. Here we have used subscripts to distinguish between the pairings over the cotangent bundle $T^*Q$ and the Lie algebra $\mathfrak{g}$. One can similarly define the $\diamond$ operator for the cotangent bundle $T^*V$. We will drop the subscripts in subsequent derivations when the space corresponding to the pairing is evident from the context. 

Upon applying the constrained variations in \eqref{eq:metamorphsis constrained variations}, the variational principle in \eqref{eq:metamorphsis var principle} takes its Euler-Poincar\'e form,
\begin{align}
    \begin{split}
        0 = \delta S &= \int \scp{\frac{\delta \ell}{\delta u}}{\delta u} + \scp{\frac{\delta \ell}{\delta n}}{\delta n} +  \scp{\frac{\delta \ell}{\delta \nu}}{\delta \nu} + \scp{\frac{\delta \ell}{\delta a}}{\delta a}\,dt\\
        &=\int \scp{\frac{\delta \ell}{\delta u}}{\p_t\eta - \ad_u\eta} + \scp{\frac{\delta \ell}{\delta n}}{w - \mathcal{L}_\eta n} + \scp{\frac{\delta\ell}{\delta \nu}}{\p_t w + \mathcal{L}_u w-\mathcal{L}_\eta \nu} + \scp{\frac{\delta\ell}{\delta a}}{-\mathcal{L}_\eta a}\,dt\\
        &=\int\scp{-\p_t\frac{\delta \ell}{\delta u} - \ad^*_u\frac{\delta \ell}{\delta u} + \frac{\delta\ell}{\delta n}\diamond n + \frac{\delta\ell}{\delta \nu}\diamond \nu + \frac{\delta\ell}{\delta a}\diamond a}{\eta} + \scp{-\p_t\frac{\delta\ell}{\delta \nu} + \mathcal{L}^T_u\frac{\delta\ell}{\delta \nu} + \frac{\delta\ell}{\delta n}}{w}\,dt \,,
    \end{split}
\end{align}
where the coadjoint operation $\ad^*: \mathfrak{g}\times \mathfrak{g}^* \to \mathfrak{g}^*$ for right action is defined by the $L^2$ pairing
\begin{align}
\scp{\ad^*_u\mu}{v} := \scp{\mu}{\ad_uv} = \scp{\mu}{-{\cal L}_uv}\,,
\quad\hbox{and}\quad
\ad^*_u\mu = {\cal L}_u\mu
\quad\hbox{for}\quad
\mu\in \mathfrak{g}^*, \quad u,v \in \mathfrak{g}
\,.
\label{def: ad-star}
\end{align}
The stationary conditions resulting from the variations, together with the definitions of $w$ and $a$, provide the evolution equations for the dynamics of the whole system,\footnote{As discussed further below, the equation set in \eqref{eq:metamorphsis EP eq} for WMFI dynamics taking place in the frame of the fluid motion closely tracks the equations for the dynamics of complex fluids reviewed authoritatively in \cite{GBR09}.}
\begin{align}
\begin{split}
    (\p_t + \ad^*_u)\frac{\delta \ell}{\delta u} &= \frac{\delta\ell}{\delta n}\diamond n + \frac{\delta\ell}{\delta \nu}\diamond \nu + \frac{\delta\ell}{\delta a}\diamond a\,, \\
    (\p_t + \mathcal{L}_u)\frac{\delta\ell}{\delta \nu} &= \frac{\delta\ell}{\delta n} \,,\\
    (\p_t + \mathcal{L}_u) n &= \nu \,,\\
    (\p_t + \mathcal{L}_u) a &= 0\,,
\end{split} \label{eq:metamorphsis EP eq}
\end{align}
where we have used the fact that $-\mathcal{L}_u^T = \mathcal{L}_u$ under integration by parts in the $L^2$ pairing. We shall refer the equations \eqref{eq:metamorphsis EP eq} as Euler-Poincar\'e equations with cocyles, versions of which have also been derived in a variety of places elsewhere for hybrid dynamics, as well as when using \emph{metamorphsis reduction} in \cite{GBHR12}. Note that the second and third equations in \eqref{eq:metamorphsis EP eq} are the Euler-Lagrange equations in the frame of reference moving with the dynamics on $\mathfrak{g}$. Hence, the usual time derivative found in the Euler-Lagrange equations has been replaced by the advective derivative $\p_t + \mathcal{L}_u$. It should also be noted that the third equation in \eqref{eq:metamorphsis EP eq} takes the same form as the kinematic boundary condition, commonly found in free boundary fluid dynamics models. Thus, the kinematic boundary constraint may be interpreted as a relationship between position and velocity in a moving frame of reference, in agreement with the statement that a particle initially on the surface remains so. See, e.g., \cite{Proceedings2022}.
\begin{remark}[Hamilton--Pontryagin principle and semidirect product reduction]
The Hamilton--Pontryagin principle equivalent to the constrained variational principle \eqref{eq:metamorphsis var principle} is the following,
\begin{align}
    0 = \delta\int \ell(u, q g^{-1}, v g^{-1}, a) + \scp{m}{\dot{g} g^{-1} - u}  + \scp{b}{a_0 g^{-1} -a} + \scp{p g^{-1}}{\Dot{q} g^{-1} - v g^{-1}}\,dt\,,
    \label{eq:metamorphsis HP var principle}
\end{align}
where all variations are arbitrary, modulo vanishing at the end points in time. Note that the kinematic constraint $\scp{p}{\dot{q} - v}$ has been acted on from the right by $g^{-1}$ and it takes the form of the kinematic boundary condition for a free boundary. Together with the constraint $\scp{m}{\dot{g}g^{-1} - u}$, one can view the tuple $(g, q)$ are elements of the semi-direct product group $S = G\circledS Q$ since the relation
\begin{align}
    \p_t(g, q)(g, q)^{-1} = (\dot{g}g^{-1}, \dot{q}g^{-1})\,,
\end{align}
is isomorphic to the Lie algebra $\tilde{\mathfrak{s}}$ of $\tilde{S}$. See, e.g., \cite{CH13} for an another application of this relation. The metamorphosis Hamilton--Pontryagin variational principle in \eqref{eq:metamorphsis HP var principle} becomes
\begin{align}
    0 = \delta\int \ell(u, n, \nu, a) + \scp{(m,\pi)}{\p_t(g, q)(g, q)^{-1}  - (u, \nu)}_{\tilde{\mathfrak{s}}} + \scp{b}{a_0 g^{-1} -a}\,dt\,, 
\end{align}
when the reduced definitions $u:=\dot{g}g^{-1}$, $n = q g^{-1}$, $\nu = \dot{q}g^{-1}$ are used, and one defines $\pi := pg^{-1}$. The subscript $\tilde{\mathfrak{s}}$ included in the pairing indicates that the pairing is to be taken with respect to $\tilde{\mathfrak{s}}$. 
\end{remark}
\begin{remark}[Symmetry breaking]
The explicit dependence of the Lagrangian, $\ell$, on $n= q g^{-1}$ means that the dynamics is not reduced by the entire symmetry group $\tilde{S} = G\circledS Q$ from the cotangent bundle $T^*\tilde{S}$. Instead, the reduction is only by $G$, so the dynamics takes place on the Lie-algebra $\Tilde{\mathfrak{s}} := \mathfrak{g}\circledS(T^*Q)$. Thus, the canonical two-cocyle arising from metamorphisis reduction of this type is inherited from the canonical Hamiltonian motion on $T^*Q$. 
\end{remark}

\begin{remark}[A composition of maps]\label{rmk:CoM}
As shown in \cite{S2022}, the Euler-Poincar\'e equations \eqref{eq:metamorphsis EP eq} can similarly be obtained from a Lagrangian depending on $TQ \times V^*$ in which an element of $TQ$ is defined as a composition. This feature builds on the `composition of maps' approach discussed in \cite{Proceedings2022}. The resulting Lagrangian is defined to be right invariant under the action of $G$ as
\begin{equation}\label{eqn:CoM_Lagrangian}
    L(g,\dot g , n g , (n g)^{\bs\cdot} , a_0) = \ell(\dot gg^{-1} , n , (n g)^{\bs\cdot} g^{-1} , a_t ) \,,
\end{equation}
where we have again denoted the composition of two maps by concatenation from the right. By writing the composition as a pullback, the Lie chain rule allows us to define $\nu$ as follows
\begin{equation}\label{eqn:LieChain}
    (g^*n)^{\bs \cdot}g^{-1} = g^*\big[ (\p_t + \mathcal{L}_u)n \big]g^{-1} = (\p_t + \mathcal{L}_u)n =: \nu \,,
\end{equation}
since the pull-back by $g$ is the inverse of the push-forward by $g$.
Indeed, we see that this agrees with the definition made in the reduction by stages process above; namely, $\nu = \dot q g^{-1}$.
\end{remark}

\subsection{The Hamiltonian formulation.}\label{subsec:Ham} One may also consider the reduced variational principle from the perspective of Hamiltonian mechanics. Indeed, the corresponding reduced Hamiltonian
\begin{align}
    h:\mathfrak{g}^*\times T^*Q \times V^* \rightarrow \mathbb{R}\,,
\end{align}
can be derived equivalently by reduction by symmetry on the Hamiltonian side. Please note that the Hamiltonian $H:T^*(G\times Q)\times V^*$ is invariant under the right action of $G$, where the group action is denoted by concatenation. The reduced Hamiltonian $h$ can be found by the quotient map
\begin{align}
    T^*(G\times Q)\times V^* \rightarrow \mathfrak{g}^*\times T^*Q\times V^*\,,\quad (g, \alpha, q, p, a_0)\rightarrow (m, n, \pi, a)\,,
\end{align}
where $m := \alpha g^{-1}$ and $\pi := p g^{-1}$. One can equivalently use the reduced Legendre transform 
\begin{align}
    h(m, n, \pi, a) = \scp{m}{u} + \scp{\pi}{\nu} - \ell(u,n,\nu,a)\,,
\end{align}
to obtain the reduced Hamiltonian $h$ from the corresponding reduced Lagrangian $\ell$. Noting that $\frac{\delta\ell}{\delta \nu} = \pi$ and $\frac{\delta h}{\delta \pi} = \nu$, one can write \eqref{eq:metamorphsis EP eq} in Hamiltonian form as
\begin{align}
\begin{split}
    \left(\p_t + \ad^*_u\right)m &=  - \frac{\delta h}{\delta n}\diamond n - \frac{\delta h}{\delta \pi}\diamond \pi - \frac{\delta h}{\delta a}\diamond a\,, \\
    \left(\p_t + \mathcal{L}_u \right)\pi &= -\frac{\delta h}{\delta n}\,,\\
    \left(\p_t + \mathcal{L}_u \right)n &= \frac{\delta h}{\delta \pi}\,,\\
    \left(\p_t + \mathcal{L}_u \right)a &= 0\,, \quad\hbox{where}\quad u:=\frac{\delta h}{\delta m}\,,
\end{split} \label{eq:metamorphsis LP eq}
\end{align}
which are the Lie-Poisson equations with cocycles. In particular, the second and third equations in \eqref{eq:metamorphsis LP eq} are \emph{Hamilton's canonical equations}, boosted into a moving frame of reference. At the level of the equations, this is equivalent to replacing the time derivative with $\p_t + \mathcal{L}_u$, as we saw with the Euler-Lagrange equations in \eqref{eq:metamorphsis EP eq}.
Hence, one can arrange \eqref{eq:metamorphsis LP eq} into Poisson bracket form as 
\begin{align}
    \p_t\begin{pmatrix}m\\a\\ \pi \\ n \end{pmatrix} =
    -\begin{pmatrix}\ad^*_{\fbox{}}m & \fbox{}\diamond a & \fbox{}\diamond \pi & \fbox{}\diamond n \\
    \mathcal{L}_{\fbox{}}a & 0 & 0 & 0\\
    \mathcal{L}_{\fbox{}}\pi & 0 & 0 & 1\\
    \mathcal{L}_{\fbox{}}n & 0 & -1 & 0\end{pmatrix}
    \begin{pmatrix}\frac{\delta h}{\delta m} = u\\ \frac{\delta h}{\delta a} = -\frac{\delta \ell}{\delta a} \\ \frac{\delta h}{\delta \pi} = \nu \\ \frac{\delta h}{\delta n} = -\frac{\delta \ell}{\delta n}\end{pmatrix}\,.\label{eq:metamorphosis Poisson matrix}
\end{align}

The Hamiltonian structure of the Poisson bracket \eqref{eq:metamorphosis Poisson matrix} is \emph{tangled} in the sense that the Lie-Poisson bracket on $\mathfrak{g}^*\circledS V^*$ is coupled to the canonical Poisson bracket on $T^*Q$ via the semidirect product structure. The Poisson structure is then $\mathfrak{g}^*\circledS V^* \circledS T^*Q$. One can \emph{untangle} the Hamiltonian structure of the Poisson bracket \eqref{eq:metamorphosis Poisson matrix} into the direct sum of the Lie-Poisson bracket on $\mathfrak{g}^*\circledS V^*$ and the canonical Poisson bracket on $T^*Q$. This is done via the map
\begin{align}
    (m, n, \pi, a) \in \mathfrak{g}^*\times T^*Q\times V^* \rightarrow (m + \pi\diamond n, n, \pi, a) =: (\kappa, n, \pi, a) \in \mathfrak{g}^*\times T^*Q\times V^*\,. \label{eq:untangling map}
\end{align}
The untangled Poisson structure can be directly calculated and written in terms of the transformed Hamiltonian $h_{HP}(\kappa, n, \pi, a)$ as 
\begin{align}
    \p_t\begin{pmatrix}\kappa\\a\\ \pi \\ n \end{pmatrix} =
    -\begin{pmatrix}\ad^*_{\fbox{}}\kappa & \fbox{}\diamond a & 0 & 0\\
    \mathcal{L}_{\fbox{}}a & 0 & 0 & 0\\
    0 & 0 & 0 & 1\\
    0 & 0 & -1 & 0\end{pmatrix}
    \begin{pmatrix}\frac{\delta h_{HP}}{\delta \kappa} = u\\ \frac{\delta h_{HP}}{\delta a} = -\frac{\delta \ell_{HP}}{\delta a} \\ \frac{\delta h_{HP}}{\delta \pi} = \nu - \mathcal{L}_u n \\ \frac{\delta h_{HP}}{\delta n} = -\frac{\delta \ell}{\delta n} + \mathcal{L}_u \pi \end{pmatrix}\,.\label{eq:untangled Poisson matrix}
\end{align}
As pointed out in \cite{GBHR12}, the untangled Poisson structure can be derived via the \emph{Hamilton--Poincar\'e} reduction principle when the Hamiltonian collectivises into the momentum map $\kappa = m + \pi\diamond n$. The dual map of \eqref{eq:untangling map} is
\begin{align}
    (u, n ,\nu, a)\in \mathfrak{g}\times TQ\times V^* \rightarrow (u, n, \nu - \mathcal{L}_u n, a) =:(u, n, \dot{n}, a) \in \mathfrak{g}\times TQ\times V^*\,,
\end{align}
which are the variables in \emph{Lagrange Poincar\'e reduction} of $L$ to the reduced Lagrangian $\ell_{LP}$ and we have the equivalence
\begin{align}
    \ell(u, n ,\nu, a) = \ell_{LP}(u, n ,\nu - \mathcal{L}_u n, a)\,.
\end{align}
\begin{remark}[Untangling from constrained variations]
    Recall the constrained variations \eqref{eq:metamorphsis constrained variations}. The choice of whether to define the variations in terms of $(\delta q) g^{-1}$ or $\delta (q g^{-1})$ will lead respectively to the tangled and untangled Euler-Poincar\'e equations corresponding to the Poisson operators \eqref{eq:metamorphosis Poisson matrix} and \eqref{eq:untangled Poisson matrix}.  This is due to the correspondence between the variations and definitions of $\nu$ and $\dot{n}$ as the tangled and untangled velocities in $TQ$.
\end{remark}

%\paragraph{Significance of the Kelvin-Noether theorem}
%We note that the untangled momentum map $\kappa = m + \pi\diamond n$ produces the $1$-form density comprising the sum of two terms in the same form as those seen in the integrand of the Kelvin-Noether theorem equation for GLM in \fbox{\eqref{GLM-circ-Lag}}.
%\todo[inline]{RH: What is the boxed equation reference? 
%\\DH: It was for a GLM equation that is no longer in the paper.}
By assuming further that the fluid density $D$ is also advected by the flow, i.e. $\p_t D + \mathcal{L}_u D = 0$, we find the following Kelvin-circulation theorem for the momentum map $\kappa = m + \pi\diamond n$,
\begin{align}
    \frac{d}{dt}\oint_{c(u)}\underbrace{ \frac{m}{D} +\frac{\pi\diamond n}{D}}_{\hbox{`Momentum shift'}}
    = -\oint_{c(u)} \frac{1}{D}\left(\frac{\delta h}{\delta a}\diamond a + \frac{\delta h}{\delta D}\diamond D\right)\,.
\label{eq:KC untangled}
\end{align}
The additional term $(\pi\diamond n)/D$ in the integrand of the Kelvin-Noether theorem in \eqref{eq:KC untangled} is a shift in momentum $1$-form, as observed earlier in the GLM and CL cases. The canonically conjugate pair $(\pi, n)$ here are  Hamiltonian variables whose dynamics takes place in the frame of the fluid motion, appearing in the result of Hamilton's principle in equation \eqref{eq:untangled Poisson matrix}. Using the tangled form of the Poisson matrix \eqref{eq:metamorphosis Poisson matrix} and the untangled Kelvin-Noether theorem \eqref{eq:KC untangled} yields the separated Kelvin-Noether equations,
\begin{align}
\begin{split}    
    \frac{d}{dt}\oint_{c(u)} \frac{m}{D} &= -\oint_{c(u)} \frac{1}{D}\left(\frac{\delta h}{\delta a}\diamond a + \frac{\delta h}{\delta D}\diamond D\right)
    - \oint_{c(u)}\underbrace{\frac{1}{D}\left( \frac{\delta h}{\delta n}\diamond n - \pi \diamond \frac{\delta h}{\delta \pi} \right)}_{\hbox{Non-inertial force}} \,,\\ 
    \frac{d}{dt}\oint_{c(u)}\frac{\pi\diamond n}{D} &= \oint_{c(u)}\frac{1}{D}\left( \frac{\delta h}{\delta n}\diamond n - \pi \diamond \frac{\delta h}{\delta \pi} \right) \,.
\end{split}
\label{eq:KC tangled}
\end{align}
Thus, the wave degree of freedom introduces a non-inertial force reminiscent of the Coriolis force, except that it has become dynamical.
Equations \eqref{eq:KC tangled} are interpreted as the result of shifting the Hamiltonian $(\pi, n)$ dynamics into the frame of the moving fluid. In the inertial Eulerian frame, the result of the Galilean shift of the Hamiltonian $(\pi, n)$ dynamics is represented by the shift in the momentum $1$-form in the Kelvin circulation integrand in \eqref{eq:KC untangled}. In the non-inertial Lagrangian frame, the result of the Galilean shift of the Hamiltonian $(\pi, n)$ dynamics is represented as the additional non-inertial force on the current in \eqref{eq:KC tangled}. %Whenever the order parameter $q$ in this shift has a spatial gradient, the shift in the momentum per unit mass will generate a force on the fluid current, just as in the case of rotation, or CL, or GLM.

\begin{remark}[Partial Legendre transform (Routhian)]
    One can show the Hamilton--Pontryagin principle in \eqref{eq:metamorphsis HP var principle} takes a form similar to that introduced in \cite{H2019} through a partial Legendre transform of a particular form of the reduced Lagrangian $\ell$. Namely, one assumes that $\ell$ is separable between the variables in $TQ$ and variables in $\mathfrak{g}\times V^*$,
    \begin{align}
        \ell(u, n ,\nu,a) = \ell_{\mathfrak{g}\times V^*}(u, a) + \ell_{TQ}(n, \nu)\,.
    \end{align}
    After using the partial Legendre transform to obtain the Hamiltonian 
    \begin{align}
        h_{T^*Q}(\pi, n) := \scp{\pi}{\nu} - \ell_{TQ}(n,\nu) \,,
    \end{align}
    one inserts $h_{T^*Q}$ into the Hamilton-Pontryagin form \eqref{eq:metamorphsis HP var principle} to find the equivalent action principle
    \begin{align}
        0 = \delta\int \ell_{\mathfrak{g}\times V^*}(u, a) + \scp{m}{\dot{g} g^{-1} - u} + \scp{b}{a_0 g^{-1} -a} + \scp{\pi}{\dot{q}g^{-1}} - h_{T^*Q}(\pi, qg^{-1}) \,dt\,. \label{eq:metamorphsis HP var partial Leg}
    \end{align}
    In terms of the $( \pi,n)$ variables, one can cast \eqref{eq:metamorphsis HP var partial Leg} into a familiar form for wave dynamics seen, e.g. in \cite{H2019}. Namely, the Hamilton-Pontryagin form \eqref{eq:metamorphsis HP var partial Leg}  can be cast as a phase-space Lagrangian,
    \begin{align}
        0 = \delta\int \ell_{\mathfrak{g}\times V^*}(u, a) + \scp{\pi}{\p_t n + \mathcal{L}_u n} - h_{T^*Q}(\pi, n) \,dt\,, 
    \label{eq:metamorphsis HP var partial Ham}
    \end{align}
    where we have introduced the constrained variations $\delta u = \p_t \eta - \ad_u \eta$ and $\delta a = - \mathcal{L}_u a$ in place of the Hamilton--Pontryagin constraints and the canonical Hamiltonian variables $( \pi,n)$ can be varied arbitrarily.
    
    Equivalently, the metamorphosis phase-space form in \eqref{eq:metamorphsis HP var partial Ham} can be seen from the perspective of the `composition of maps' form of the Lagrangian discussed in Remark \ref{rmk:CoM}. Indeed, beginning from the Lagrangian \eqref{eqn:CoM_Lagrangian}, notice that the form of the right hand side of the inner product term of equation \eqref{eq:metamorphsis HP var partial Ham} is a direct consequence of equation \eqref{eqn:LieChain}.
    For completeness, we compute the stationary condition of \eqref{eq:metamorphsis HP var partial Ham} which gives the equations of motion
    \begin{align}
    \begin{split}
        (\p_t + \ad^*_u)\frac{\delta \ell_{\mathfrak{g}\times V^*}}{\delta u} &= - \frac{\delta h_{T^*Q}}{\delta n}\diamond n - \frac{\delta h_{T^*Q}}{\delta \pi}\diamond \pi + \frac{\delta\ell_{\mathfrak{g}\times V^*}}{\delta a}\diamond a\,, \\
        (\p_t + \mathcal{L}_u)\pi &= -\frac{\delta h_{T^*Q}}{\delta n} \,,\\
        (\p_t + \mathcal{L}_u) n &= \frac{\delta h_{T^*Q}}{\delta \pi} \,,\\
        (\p_t + \mathcal{L}_u) a &= 0\,, \label{eq:metamorphsis LP partial Leg}
    \end{split}
    \end{align}
    which are equivalent to the general Lie-Poisson equations with cocycles \eqref{eq:metamorphsis LP eq}. Using the evolution equations of $\pi$ and $n$ in \eqref{eq:metamorphsis LP partial Leg}, we can calculated explicitly the evolution of $h_{T^*Q}$ to be
    \begin{align}
    \begin{split}
        \frac{d}{dt}h_{T*Q}(\pi, n) &= \scp{\p_t \pi + \mathcal{L}_u \pi}{\frac{\delta h_{T*Q}}{\delta \pi}} + \scp{\p_t n + \mathcal{L}_u \pi}{\frac{\delta h_{T*Q}}{\delta n}} \\
        &= \scp{-\frac{\delta h_{T*Q}}{\delta n}}{\frac{\delta h_{T*Q}}{\delta \pi}} + \scp{\frac{\delta h_{T*Q}}{\delta \pi}}{\frac{\delta h_{T*Q}}{\delta n}} = 0\,,
    \end{split}
    \end{align}
    which is the conservation of $h_{T^*Q}$ along the flow generated by $u$.
\end{remark}
% \begin{remark}[Affine action construction of two cocycle]
%     To obtain the two cocycle in \eqref{eq:metamorphosis Poisson matrix}, the cotangent bundle structure of $(p,q)$ is assumed. However, the cotangent bundle assumption can be interchanged with affine actions which recovers the same Poisson structure. The construction 
% \end{remark}

\subsection{Additional symmetry}\label{subsec:additional_symmetry} So far, we have only considered the case where the symmetries of the system exist solely in the Lie group $G$. It is natural to extend the reduction principle to consider cases where the configuration manifold $Q$ is also a Lie group with corresponding Lie algebra $\mathfrak{q}$. Let $G$ act on $Q$ as a Lie group homomorphism, i.e., $G\times Q\to Q$. Additionally, we introduce explicit dependence of an order parameter $\chi_0 \in V_Q^*$ for $Q$ to the Lagrangian $L$ such that
\begin{align}
    L(g,\dot{g}, q,\dot{q}, \chi_0, a_0) : TG\times TQ \times V^*_Q\times V^*\rightarrow \mathbb{R}\,, \label{eq:unreduced lag additional symm}
\end{align}
and assume that the Lagrangian is invariant under the action of both $Q$ and $G$. For simplicity of exposition, let us consider only the right action of $q\in Q$ on $TQ$ and $\chi_0$; so the $Q$-reduced Lagrangian, $\Tilde{L}$, takes the following form
\begin{align}
    L(g, \dot{g}, q, \dot{q}, \chi_0, a_0) =: \Tilde{L}(g, \dot{g}, \dot{q}q^{-1}, \chi_0 q^{-1}, a_0): TG\times\mathfrak{q}\times V_Q^*\times V^* \rightarrow \mathbb{R}\,.
\end{align}
After the reduction by $Q$, the equations of motions are Lagrange-Poincar\'e equations \cite{CMR2001}. The further reduction by $G$ then defines the fully reduced Lagrangian $\tilde{\ell}$ by
\begin{align}
    \tilde{L}(g, \dot{g}, \dot{q}q^{-1},\chi_0 q^{-1}, a_0) 
    &= \tilde{L}(e, \dot{g}g^{-1}, (\dot{q}q^{-1})g^{-1},(\chi_0 q^{-1})g^{-1}, a_0g^{-1})
    \\&=: \tilde{\ell}(u, \omega, \chi, a): \mathfrak{g}\times\mathfrak{q}\times V_Q^* \times V^* \rightarrow \mathbb{R}\,, \label{eq:Q G reduced Lag}
\end{align}
where one defines the following abbreviated notation,
\begin{align}
    u:=\dot{g}g^{-1}, \quad 
    \omega := (\dot{q}q^{-1})g^{-1}, \quad
    \chi := (\chi_0 q^{-1})g^{-1}, \quad
    \hbox{and}\quad
    a := a_0g^{-1}\,.
\label{EulVar-defs}
\end{align}
The reduced Euler-Poincar\'e variational principle becomes
\begin{align}
    0 = \delta S = \delta \int_{t_0}^{t_1} \tilde{\ell}(u, \omega, \chi, a)\,dt\,, \label{eq:Q G reduced var prin}
\end{align}
subject the constrained variations obtained from the definitions of $u$, $\omega$ and $a$ in equation \eqref{EulVar-defs},
\begin{align}
\begin{split}
    \delta u &= \p_t\eta - \ad_u \eta \,,\\
    \delta \omega &= \p_t u - \mathcal{L}_\eta \omega + \mathcal{L}_u \gamma - \ad_\omega \gamma\,,\\
    \delta \chi &= -\mathcal{L}_\eta \chi - \widehat{\mathcal{L}}_\gamma \chi\,,\\
    \delta a &= -\mathcal{L}_\eta a\,.
\end{split}\label{eq:Q G reduced variations}
\end{align}
Here, we denote $\gamma := (\delta q q^{-1})g^{-1}$ and $\eta$ is chosen arbitrarily and vanishes at the endpoints $t=t_0,t_1$. We also introduce the notation $\widehat{\mathcal{L}}_{\gamma}$ for the action of an arbitrary Lie algebra element $\gamma \in \mathfrak{q}$. As in the definition of the diamond operator $(\diamond)$ in \eqref{def-diamond} for the Lie-derivative of vector fields $\eta\in \mathfrak{g}$, we define the diamond operator $(\widehat{\diamond})$ with respect to the action by vector fields $\mathfrak{q}$ through
\begin{align}
    \scp{-\xi \,\widehat{\diamond}\,\chi}{\gamma}_{\mathfrak{q}} = \scp{\xi}{\widehat{\mathcal{L}}_\gamma \chi}_{T^*V_Q}\,,
\end{align}
for all $(\xi, \chi)\in T^*V_Q$ and $\gamma\in \mathfrak{q}$. After taking variations one finds the Euler-Poincar\'e equations from the reduced Euler-Poincar\'e principle \eqref{eq:Q G reduced var prin} as
\begin{align}
    \begin{split}
        \p_t \frac{\delta \tilde{\ell}}{\delta u} + \ad^*_u\frac{\delta \tilde{\ell}}{\delta u}  &= \frac{\delta \tilde{\ell}}{\delta \omega}\diamond \omega + \frac{\delta \tilde{\ell}}{\delta a}\diamond a + \frac{\delta \tilde{\ell}}{\delta \chi}\diamond \chi\,,\\
        \p_t \frac{\delta \tilde{\ell}}{\delta \omega} + \ad^*_\omega \frac{\delta \tilde{\ell}}{\delta \omega} + \mathcal{L}_u \frac{\delta \tilde{\ell}}{\delta \omega} &=  \frac{\delta \tilde{\ell}}{\delta \chi}\widehat{\diamond} \chi \,,\\
        \p_t\chi + \mathcal{L}_u \chi + \widehat{\mathcal{L}}_\omega \chi &= 0\,,\\
        \p_t a + \mathcal{L}_u a &= 0\,.
    \end{split}\label{G Q metamorphsis EP eq}
\end{align}
% Note that in \eqref{eq:Q G reduced variations} we have abused the notation of $\mathcal{L}$ to denote both the Lie-derivative of the vector fields $\eta \in \mathfrak{g}$ and the action of arbitrary Lie algebra $\gamma \in \mathfrak{q}$. Based on the first argument of $\mathcal{L}_{\fbox{}}$, the meaning of notation should be self evident from the context. Similarly for the diamond operator $\diamond$, the image of $\diamond$ is either $\mathfrak{g}^*$ or $\mathfrak{q}^*$ which is evident clear from the context.
Under similar considerations on the Hamiltonian side, we can construct the reduced Hamiltonian $\tilde{h}(m,\lambda ,a):\mathfrak{g}^*\times\mathfrak{q}^*\times V_Q^*\times V^*\rightarrow\mathbb{R}$ via the Legendre transform such that $\lambda := \frac{\delta \tilde{\ell}}{\delta \omega}$ and $m := \frac{\delta \tilde{\ell}}{\delta u}$. The equations \eqref{G Q metamorphsis EP eq} can then be written in a Poisson matrix form
% \begin{align}
% \begin{split}
%     &\p_t m + \ad^*_u m + \lambda\diamond \frac{\delta \tilde{h}}{\delta \lambda} + \frac{\delta \tilde{h}}{\delta a}\diamond a = 0\,, \\
%     &\p_t\lambda + \ad^*_\omega \lambda - \mathcal{L}^T_u \lambda = 0\,,
% \end{split}
% \end{align}
\begin{align}
    \p_t\begin{pmatrix}m\\ a \\\lambda  \\ \chi \end{pmatrix} =
    -\begin{pmatrix}\ad^*_{\fbox{}}m & \fbox{}\diamond a & \fbox{}\diamond \lambda & \fbox{}\diamond \chi\\
    \mathcal{L}_{\fbox{}}a & 0 & 0 & 0\\
    \mathcal{L}_{\fbox{}}\lambda & 0& \ad^*_{\fbox{}}\lambda & \fbox{}\,\widehat{\diamond}\,\chi \\
    \mathcal{L}_{\fbox{}}\chi & 0 & \widehat{\mathcal{L}}_{\fbox{}}\chi & 0
    \end{pmatrix}
    \begin{pmatrix}\frac{\delta \tilde{h}}{\delta m} = u\\ \frac{\delta \tilde{h}}{\delta a} = -\frac{\delta \tilde{\ell}}{\delta a}\\ \frac{\delta \tilde{h}}{\delta \lambda} = \omega \\ \frac{\delta \tilde{h}}{\delta \chi} = -\frac{\delta \tilde{\ell}}{\delta \chi}\end{pmatrix}\,. \label{G Q metamorphsis PB}
\end{align}
The Lie-Poisson matrix in equation \eqref{G Q metamorphsis PB} defines a Lie-Poisson bracket on $\mathfrak{g}^*\times\mathfrak{q}^*\times V_Q^* \times V^*$, which is the same as the bracket on the dual of the  semidirect product Lie algebra $\mathfrak{s}^* = \mathfrak{g}^*\circledS((\mathfrak{q}^*\circledS V_Q^*)\oplus V^*)$.
% that is 
% \begin{align}
%     \{f,g\}_{\mathfrak{s}^*} = -\scp{m}{\left[\frac{\delta f}{\delta m}, \frac{\delta g}{\delta m}\right]} + \scp{a}{\frac{\delta f}{\delta a}\frac{\delta g}{\delta m} - \frac{\delta g}{\delta a}\frac{\delta f}{\delta m}} - \scp{\lambda}{\left[\frac{\delta f}{\delta \lambda}, \frac{\delta g}{\delta \lambda}\right] + \frac{\delta f}{\delta \lambda}\frac{\delta g}{\delta m} - \frac{\delta g}{\delta \lambda}\frac{\delta f}{\delta m} }\,.
% \end{align}
Thus, equations \eqref{G Q metamorphsis PB} are the canonical Lie-Poisson equations on $\mathfrak{s}$, the Lie-algebra of the semi-direct product group $S = G\circledS((Q\circledS V_Q)\oplus V)$, under the reduction by symmetry of $S$ itself. \\ \\ 
Reduction by left action follows an analogous procedure, and a combination of left and right reduction may also be applied. An extensive literature exists for reduction by symmetry in the theory and applications of geometric mechanics, whose foundations are reviewed in Abraham and Marsden \cite{AM78}.
\begin{remark}
    A geophysical fluid system with similar Poisson structure to \eqref{G Q metamorphsis PB} arises in the vertical slice models \cite{CH13}. This model is of the form corresponding to equation \eqref{eq:unreduced lag additional symm}, where $G = \Diff(\mathcal{D})$ for $\mathcal{D}\in\mathbb{R}^2$, $Q = \Diff(\mathcal{D})$, and there are no symmetry breaking order parameters $\chi_0$. Then, the reduction process gives $\omega \in \mathfrak{X}$ and $\pi\in \mathfrak{X}^*$ and the Lie-Poisson matrix becomes, 
\begin{align}
    \p_t\begin{pmatrix}m\\ \pi\\a \end{pmatrix} =
    -\begin{pmatrix}\ad^*_{\fbox{}}m & \ad^*_{\fbox{}}\pi& \fbox{}\diamond a  \\
    \ad^*_{\fbox{}}\pi & \ad^*_{\fbox{}}\pi & 0 \\
    \mathcal{L}_{\fbox{}}a & 0 & 0 
    \end{pmatrix}
    \begin{pmatrix}\frac{\delta h}{\delta m} = u\\ \frac{\delta h}{\delta \pi} = \omega \\ \frac{\delta h}{\delta a} = -\frac{\delta l}{\delta a}\end{pmatrix}\,.
\end{align}
\end{remark}

Let us summarise the two types of reduction procedure explained in this section. Starting from a Lagrangian $L$ defined on $T(G\times Q)\times V^*$ the two reduction pathways discussed here can be represented diagrammatically as 
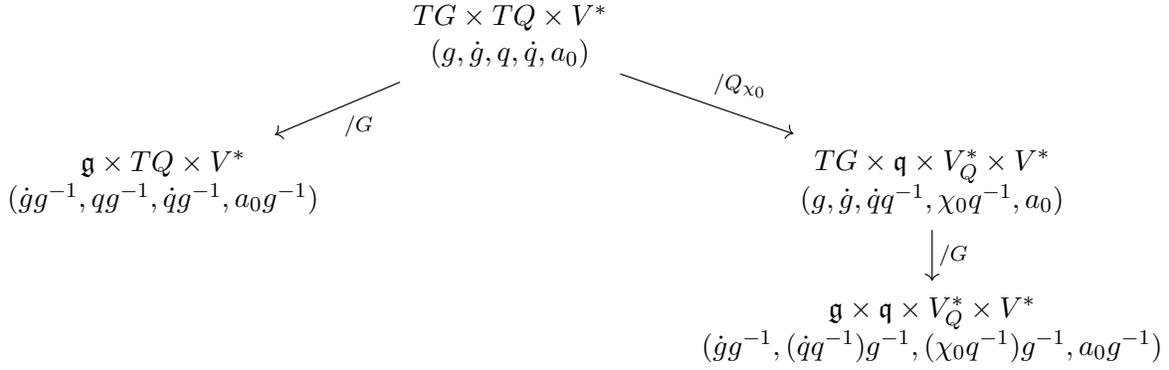
\begin{figure}[H]
\centering
\begin{tikzcd}
&
\begin{matrix}
TG \times TQ\times V^*\\
(g,\dot{g},q,\dot{q},a_0)
\end{matrix}
\arrow[dr, "/Q_{\chi_0}"]
\arrow[dl, "/G"]
&
\\
\begin{matrix}
\mathfrak{g} \times TQ\times V^*\\
(\dot{g}g^{-1},qg^{-1},\dot{q}g^{-1},a_0g^{-1})%=: (u, n, \nu, a)
\end{matrix}
&
&
\begin{matrix}
TG \times \mathfrak{q}\times V_Q^*\times V^*\\
(g,\dot{g},\dot{q}q^{-1},\chi_0 q^{-1},a_0)
\end{matrix}
\arrow[d, "/G"]
\\
&
&
\begin{matrix}
\mathfrak{g} \times \mathfrak{q}\times V_Q^*\times V^*\\
(\dot{g}g^{-1},(\dot{q}q^{-1})g^{-1},(\chi_0 q^{-1})g^{-1},a_0g^{-1})% =: (u, \omega, \chi, a)
\end{matrix}
\end{tikzcd}
\caption{Reduction pathways: The leftmost pathway in the diagram above describes the result of reducing the dynamics in $G$ by symmetry whilst $G$ acts on $Q$ via homomorphism. The rightmost pathway summarises the case where additional symmetry exists in $Q$, which is then assumed to also be a Lie group. In this case the Lagrangian is invariant under the action of the isotropy group, $Q_{\chi_0}$, of the additional parameter $\chi_0$ introduced in Section \ref{subsec:additional_symmetry}.}
\label{fig:reduction path}
\end{figure}
Both branches of this diagram reflects the reduction process relating to the specific WMFI models discussed in Section \ref{sec:3-NLS}. 

\section{Examples: Eulerian wave elevation field equations}\label{sec:3-NLS}

This section provides worked examples of wave mean flow interaction (WMFI) models. To better understand the structure the forthcoming models, see also Appendix \ref{appendix:SHO}, where one can find an elementary example demonstrating the coupling of a field of simple harmonic oscillators to an Euler fluid. 

\subsection{WKB internal waves in the Euler--Boussinesq (EB) approximation}
The Generalised Lagrangian Mean (GLM) theory of Andrews and McIntyre \cite{AM1978} can be expressed through the decomposition the Lagrangian trajectory of fluid particles into the mean and fluctuating part in the following way. Let the full fluid trajectory on $\mathcal{M}\subset\mathbb{R}^3$ be expressed by $\mb{X}_t(\bx_0)=g_t\mb{x}_0$ for initial position $\mb{x}_0$ and $g_t \in \Diff(\mathcal{M})$, which in turn is defined as the composition, $g_t=\Xi_t\circ\bar{g}_t$ of two diffeomorphism $\bar{g}_t$ and $\Xi_t$. That is, the Lagrangian trajectory can be written in coordinates as
\begin{align}
    \mb{X}_t =g_t\bx_0= \Xi_t \circ \bar{g}_t\bx_0 =: \bx_t +\alpha\bs{\xi}(\bx_t,t) \,, \quad\text{for}\quad \mb{x}_t := \bar{g}_t \bx_0 \,,\quad \text{and}\quad \alpha \ll 1\,, \label{eq:lag path decomp}
\end{align}
where we have chosen a particular form of $\Xi_t$ when \emph{acting on coordinates} to be $\Xi_t \bx := \bx + \alpha \bs{\xi}_t(\bx)$ for some $\bs{\xi}_t \in \mathcal{F}^3(\mathcal{M})$ with zero mean. 
The slow (mean flow) and fast (fluctuation) parts of the fluid flow are then clearly represented in the above formulae to be $\bx_t := \bar{g}_t \bx_0$ and $\alpha\bs{\xi}_t(\bx_t)$ respectively where $\bs{\xi}_t$ is interpreted as the fluctuation displacement away from the mean flow. Furthermore, the dimensionless parameter $0 < \alpha \ll 1$ is relative scale of the fluctuation dynamics against the mean flow.
Taking time derivative of the decomposed Lagrangian trajectory \eqref{eq:lag path decomp} yields the total velocity decomposition
\begin{align}
\begin{split}
    \bU_t(\mb{X}_t) := \frac{d\mb{X}_t}{dt}=\dot{g}_t g_t^{-1}\mb{X}_t 
    &= \frac{d\bx_t}{dt} + \alpha \Big(\p_t \bs{\xi}(\bx_t,t)+\frac{\p\bs{\xi}}{\p x^j_t}\frac{dx_t^j}{dt} \Big) 
    \\&= \dot{\bar{g}}_t\bar{g}_t^{-1}\bx_t + \alpha \Big(\p_t \bs{\xi}(\bx_t,t)+\frac{\p\bs{\xi}}{\p x^j_t}\cdot (\dot{\bar{g}}_t\bar{g}_t^{-1}x_t^j) \Big)  \\ 
    % \Big(\hbox{For}\quad \bU_t:=\dot{g}_t g_t^{-1}\mb{X}_t\,,\quad \bu^L:=\dot{\bar{g}}_t\bar{g}_t^{-1}\bx_t \Big) \
    &=: \bu^L(\bx_t,t) + \alpha \left( \p_t \bs{\xi} (\bx_t,t) + \bu^L\cdot\nabla_{\bx_t} \bs{\xi}(\bx_t,t) \right)
    \\& =: \bu^L(\bx_t,t) + \alpha \frac{d}{dt} \bs{\xi}(\bx_t,t) 
    \,,
\end{split} \label{eq:LagVel-decomp}
\end{align}
which is the sum of the Lagrangian mean velocity $\bu^L :=\dot{\bar{g}}_t\bar{g}_t^{-1}\bx_t$ and velocity fluctuations $\alpha\frac{d}{dt}\bs{\xi}(\bx_t, t)$. The Eulerian velocity decomposition \eqref{eq:LagVel-decomp} is the same as the GLM velocity decomposition in \cite{AM1978}
\begin{align*}
   \bU(\bx+\alpha\bs{\xi}(\mb{x},t) , t)= \mb{u}_L(\mb{x},t) + \alpha\frac{d}{dt}\bs{\xi}(\mb{x},t)  \,,\quad\hbox{where}\quad
\overline{\bU(\bx+\alpha\bs{\xi}(\mb{x},t) )} =: \mb{u}_L(\mb{x},t) 
\,.
\end{align*}

Gjaja and Holm \cite{GH1996} closed the Generalised Lagrangian Mean (GLM) theory for the case that the displacement fluctuation in $\bs{\xi}(\mb{x},t)\in \mbb{R}^3$ is given by a single-frequency travelling wave with slowly varying complex vector amplitude,
\begin{align*}
\bs{\xi}(\mb{x},t) = \frac12\Big(\mb{a}(\mb{x},t)e^{i\phi(\mb{x},t)/\epsilon}
+ \mb{a}^*(\mb{x},t)e^{-i\phi(\mb{x},t)/\epsilon}\Big)
\quad\hbox{with}\quad \epsilon\ll1
\,.
\end{align*} 
Here, $\epsilon$ is a dimensionless number introduced to distinguish between the time scales for the waves and currents by modifying the phase.
Holm \cite{H2021} simplified the wave mean flow interaction (WMFI) closure in \cite{GH1996} by neglecting pressure expansion and Coriolis force in the dispersion relation, thereby placing the WMFI theory into the present hybrid formulation, by coupling Lagrangian mean EB fluid equations to leading order Hamiltonian wave dynamics in the following variational principle
\begin{align}
\begin{split}
0 &= \delta \int_{t_0}^{t_1}
\ell_{EB}(\bu^L, D, \rho; p) - \scp{N}{\p_t + \bu^L\cdot\nabla\phi} - H_W(N,\bk)\,dt \\
& = \delta \int_{t_0}^{t_1} \int_\mathcal{M}\frac{D}{2}\big| \bu^L \big|^2 + D\bu^L\cdot \bR(\bx) - gDbz - p(D-1) \\
& \hspace{3cm} - N(\partial_t\phi  + \bu^L\cdot\nabla \phi)\,d^3x\,dt - H_W(N,\bk) \,dt\, \quad\textrm{with}\quad \bk\cdot d\bx:={\rm d} \phi\,,
\end{split}
\label{eq:WKB EB Lag}
\end{align}
where the various variables appearing in \eqref{eq:WKB EB Lag} are defined as follows. $\bu^L(\bx, t)$ denotes the coefficients of the vector field $u^L := \bu^L\cdot \nabla \in \mathfrak{X}(\mathcal{M})$ which is the Lagrangian mean velocity vector field, ${\rm curl}\,\bR(\bx)=2\bs{\Omega}(\bx)$ is the prescribed Coriolis parameter, $Dd^3x \in \text{Den}(\mathcal{M})$ is the volume element, $b \in \mathcal{F}(\mathcal{M})$ is the scalar buoyancy and $p\in\mathcal{F}(\mathcal{M})$ is the Lagrange multiplier enforcing incompressibility. As for the wave variables, $Nd^3x \in \text{Den}(\mathcal{M})$ is the wave action density and $\phi \in \mathcal{F}(\mathcal{M})$ is the canonically conjugate scalar wave phase.
% Thus, the first summand governs the Lagrangian mean EB fluid dynamics and the second summand in that variational principle governs the dynamics of the leading order fluctuations away from the mean. 
Variations of the fluid variables in \eqref{eq:WKB EB Lag} are constrained, namely $\delta u^L = \p_t w - \ad_{u^L} w$, $\delta (D d^3x) = -\mathcal{L}_{w} (Dd^3x) =- {\rm div}(D \mb{w})d^3x$ and $\delta b = -\mathcal{L}_{w} = - \mb{w}\cdot\nabla b$ for arbitrary variation $\mb{w}\cdot \nabla = w \in \mathfrak{X}(\mathcal{M})$ vanishing at endpoints. The waves variables, $N$ and $\phi$, instead have arbitrary variations which vanishes at endpoints. 

Let us consider the proposed variational principle \eqref{eq:WKB EB Lag} in the context of various Lagrangian reduction techniques discussed in \autoref{sec:2-Lagrangian reduction}.
The first summand is the reduced Lagrangian mean EB fluid Lagrangian where $u^L = \dot{g}g^{-1}$, $b_t = b_0 g^{-1}$ and $Dd^3x = (D_0d^3x_0)g^{-1}$ for $g \in \Diff(\mathcal{M})$ where the action by $g$ is the push-forward defined in \eqref{advec-reln}.
The second summand gives the phase space Lagrangian of the fluctuation variables $N$ and $\phi$, which are again, the push-forward of the same wave quantities defined on the Lagrangian mass coordinates.
That is, $N_t(\bx_t) = (N_t g_t)(\bx_0)$ and $\phi_t(\bx_t) = (\phi_t g_t)(\bx_0)$. The particular form of the variational principle \eqref{eq:WKB EB Lag} is that considered in \eqref{eq:metamorphsis HP var partial Ham} where $G = \Diff(\mathcal{M})$, $V^* = \mathcal{F}(\mathcal{M})\oplus \text{Den}(\mathcal{M})$ and $Q = \mathcal{F}(\mathcal{M})$, up to a minus sign in the kinematic constraint of $\phi$ due to convention. Thus, the equations of motion derived from the stationary conditions of \eqref{eq:WKB EB Lag} can be obtained by substituting the appropriate variational derivatives into the general form \eqref{eq:metamorphsis LP partial Leg} as the following.
The modified canonical Hamilton's equations for the wave dynamics are 
\begin{align}
(\partial_t  + \mathcal{L}_{u^L})  \phi + \frac{\delta H_W}{\delta N} = 0 \,\quad\text{and}\quad
(\partial_t  + \mathcal{L}_{u^L}) (N\,d^3x) + d \left(\frac{\delta H_W}{\delta (\bk\cdot d\bx)}\right) = 0\,,\label{eq:WKB wave Ham eq}
\end{align}
where we see the fluid velocity $u^L$ transports the wave dynamics in the reference frame of the fluid flow. The equations \eqref{eq:WKB wave Ham eq} can be assembled to give the evolution equation of the wave momentum density $N\nabla\phi\cdot d\mb{x}\otimes d^3x$ as the following, 
\begin{align}
    (\partial_t  + \mathcal{L}_{u^L}) \Big( N\nabla \phi\cdot d\bx \otimes d^3x\Big) &= -\left(
{\rm div} \Big(\frac{\delta H_W}{\delta  \bk} \Big) d\phi
- N d \Big( \frac{\delta H_W }{\delta  N}\Big)\right)\otimes d^3x \,.
\label{eq:WKB wave momentum eq}
\end{align}
The evolution of the equation of the fluid advected quantites and the evolution of the total momentum can also be derived from the variational principle to be
\begin{align}
\begin{split}
&(\partial_t  + \mathcal{L}_{u^L}) \big(\bM\cdot d\bx \otimes d^3x\big) = 
\left( 
D d \pi + Dgz d b
\right)\otimes d^3x
\,,\\&
(\partial_t  + \mathcal{L}_{u^L})  (D\,d^3x)  = 0
\,,\qquad
D=1
\,,\qquad
(\partial_t  + \mathcal{L}_{u^L}) b = 0\,,
\end{split}
\label{SVP3-det-redux}
\end{align}
where the coefficients of the Eulerian total momentum density $\bM$ and pressure $\pi$ in equation \eqref{SVP3-det-redux} are given by, 
\begin{align}
\bM := D (\bu^L + \bR(\bx)) - N\nabla \phi
\,,\qquad
\pi := \frac12 |\bu^L|^2 + \bR(\bx)\cdot \bu^L - gbz - p
\,.
\end{align}
Note that the dynamics of $\mb{M}\cdot d\mb{x}$ is independent of the form of the wave Hamiltonian $H_W$, thus one finds the Kelvin circulation dynamics of $\mb{M}\cdot d\mb{x}$,
\begin{align}
\frac{d}{dt} 
\oint_{c(\mb{u}^L)} \Big(\bu^L + \bR(\bx) - \frac{N\nabla\phi}{D} \Big)\cdot d\bx
= \oint_{c(\mb{u}^L)}  (\nabla\pi + gz \nabla b) \cdot d\bx 
\,,
\label{SALT-SNWP-GLM-total}
\end{align}
where $c(\mb{u}^L)$ is a material loop moving with the flow at velocity $\bu^L(\bx,t)$. 
The total momentum density $\bM = D (\bu^L + \bR(\bx)) - N\nabla \phi$ decomposes into the \emph{sum} of the momentum densities for the two degrees of freedom, namely, the wave and fluid degrees of freedom. Defining the fluid momentum $\mb{m}\cdot d\mb{x} := \left(\mb{u}^L + \mb{R}(\mb{x})\right)\cdot d\mb{x}$, one finds its evolution as the differences of \eqref{eq:WKB wave momentum eq} and \eqref{SVP3-det-redux}
\begin{align}
(\partial_t  + \mathcal{L}_{u^L}) \big(\mb{m}\cdot d\bx \otimes d^3x\big) = 
\left( 
D d \pi + Dgz d b
\right)\otimes d^3x - \left(
{\rm div} \Big(\frac{\delta H_W}{\delta  \bk} \Big) d\phi
- N d \Big( \frac{\delta H_W }{\delta  N}\Big)\right)\otimes d^3x 
\end{align}

{\bf WKB wave Hamiltonian in 3D.} Suppose for $H_W$ one takes the WKB wave Hamiltonian in 3D, whose variational derivatives are given by familiar wave quantities, 
\begin{align}
H_W =  \int_M N \omega(\bk) \,d^3x 
\,,\quad\hbox{with}\quad
 \frac{\delta H_W }{\delta  N}\Big|_{\bk} =  \,\omega(\bk)
 \,,\quad\hbox{and}\quad
 \frac{\delta H_W}{\delta  \bk}\Big|_{N} =  N \frac{\partial \omega(\bk) }{\partial \bk} =: \,N \bv_G(\bk)
\,,\label{separatedWaveHam-redux}
\end{align}
in which $\bv_G(\bk):=\partial \omega(\bk) / \partial \bk$ is the group velocity for the dispersion relation $\omega=\omega(\bk)$ between wave frequency, $\omega$, and wave number, $\bk$.  
Then, the explicit form of the dynamics of the WKB wave momentum $\frac{N}{D}\nabla\phi\cdot d\mb{x}$ from \eqref{eq:WKB wave momentum eq} appears as
\begin{align}
    \left(\p_t + \mathcal{L}_{u^L}\right) \left(\frac{N}{D}\nabla \phi \cdot d\bx\right) = -\frac{1}{D} \bigg(\bk \,{\rm div}\Big(N\bv_G(\bk)\Big) - N \nabla  \omega(\bk)\bigg)\cdot d\bx\,,
\end{align}
where $\nabla\omega(\bk)$ refers to a \emph{spatial} gradient of the frequency. Likewise, one has the explicit form of the Kelvin-circulation dynamics for the Eulerian fluid momentum $m = \left(\mb{u}^L + \mb{R}(\mb{x})\right)\cdot d\mb{x}$ and wave momentum $\frac{N}{D}\nabla \phi\cdot d\mb{x}$ as \eqref{SVP3-det-redux}
\begin{align}
\begin{split}
\frac{d}{dt}  \oint_{c(\mb{u}^L)} \big(\bu^L + \bR(\bx)\big)\cdot d\bx 
& = \oint_{c(\mb{u}^L)} \big(\nabla \pi + gz \nabla b\big) \cdot d\bx - 
\underbrace{\
\frac{1}{D}  \bigg(\bk \,{\rm div} \Big( N \bv_G(\bk)\Big) - N \nabla \omega(\bk)\bigg) \
}_{\hbox{WKB Wave Forcing}}\hspace{-1mm}
\cdot \,d\bx \,,\\
\,\frac{d}{dt} \oint_{c(\mb{u}^L)} \frac{N}{D}\nabla \phi \cdot d\bx 
&= -\oint_{c(\mb{u}^L)} \frac{1}{D} \bigg(\bk \,{\rm div}\Big(N\bv_G(\bk)\Big) - N \nabla  \omega(\bk)\bigg)\cdot d\bx
\end{split}
\label{Det-GLM-Kelvin}
\end{align}
where $c(\mb{u}^L)$ is a material loop moving with the flow at velocity $\bu^L(\bx,t)$. 
\begin{remark}[Summary of WKB internal wave dynamics in the Euler--Boussinesq (EB) approximation]$\,$\vspace{-4mm}

$\bullet$ \quad Equations \eqref{Det-GLM-Kelvin} and \eqref{SALT-SNWP-GLM-total} provide an additive decomposition the Kelvin circulation theorem representation of WCI in the example of EB flow. This result from the variational principle for WCI dynamics in \eqref{eq:WKB EB Lag} fits well with the vast literature of mean flow interaction. See, e.g., \cite{Peregine1976,Whitham2011,Buhler2014}.  

$\bullet$ \quad The total potential vorticity (PV) is conserved on Lagrangian mean particle paths. That is,
\begin{align} 
\partial_t Q + \bu^L\cdot \nabla Q = 0\,,
\label{GLM-PV}
\end{align}
where PV is defined as $Q := D^{-1} \nabla b \cdot {\rm curl \big(\bu^L + \bR(\bx) - D^{-1}N\nabla\phi\big)}$ with $D=1$.

$\bullet$ \quad For the WKB wave Hamiltonian in \eqref{separatedWaveHam-redux}, the phase-space Lagrangian in \eqref{eq:WKB EB Lag} has produced a model of wave interactions with the mean EB fluid current in which the total circulation separates into a sum of wave and current components. 

$\bullet$ \quad In particular, the total momentum density in the model $\bM = D (\bu^L + \bR(\bx)) - N\nabla \phi$ represents the \emph{sum} of the momentum densities for the current and wave degrees of freedom, respectively. 

$\bullet$ \quad The result from the first formula in \eqref{Det-GLM-Kelvin} implies that the WKB wave contribution can feed back to create circulation of the fluid current. However, if waves are initially absent, the fluid current cannot subsequently create waves. 

$\bullet$ \quad The latter conclusion supports the interpretation of the model that the fluid variables describe mean flow properties. 
\end{remark}
The next example will consider a two-dimensional case when the wave Hamiltonian $H(N,\bk)$ corresponds to the nonlinear Schr\"odinger (NLS) equation. 

\subsection{Coupling to the nonlinear Schr\"odinger (NLS) equation}\label{subsec:NLS}

As explained in Stuart and DiPrima \cite{StuartDiprima1978}, 2D surface wave dynamics near the onset of instability may be approximated by the solutions of the NLS equation. The NLS equation is written in terms of a complex wave amplitude, $\psi$, defined in a certain Hilbert space, $\mathcal{H}$, as
\begin{align}
    i\hbar\p_t \psi = -\frac{1}{2}\Delta \psi + \kappa |\psi|^2\psi\,. 
    \label{det-NLS}
\end{align}
%\todo[inline,color=yellow]{DH: What values of $\rhbar$ and ${\color{red}m}$ are appropriate for SWOT on the sea?}
The sign of the real parameter $\kappa$ in \eqref{det-NLS} controls the behaviour of NLS solutions. In what follows, we shall use the Dirac-Frenkel (DF) variational principle pioneered in \cite{DF1934} to derive the NLS equation from Hamilton's principle and then couple its solutions to a fluid flow. The DF variational principle for the \emph{linear} Schr\"odinger equation $i\hbar\p_t\psi = \widehat{H}\psi$ with Hamiltonian operator $\widehat{H}$ can be written in the form of a \emph{phase space Lagrangian}, as
\begin{align}
    0 = \delta S = \delta \int^a_b \scp{\psi}{i\hbar\p_t\psi - \widehat{H}\psi}\,dt\,. 
    \label{eq:DF var prin}
\end{align}
The pairing $\scp{\cdot}{\cdot}$ in \eqref{eq:DF var prin}  is defined by 
\begin{align}
    \scp{\psi_1}{\psi_2} = \Re\bracket{\psi_1}{\psi_2}\,, \label{eq:complex pairing}
\end{align}
in which the bracket $\bracket{\psi_1}{\psi_2}$ is the natural inner product in Hilbert space $\mathcal{H}$. In the spatial domain $\mathcal{D}\in \mathbb{R}^2$, for the case $\mathcal{H} = L^2(\mathcal{D})$, the inner product is given by
\begin{align}\label{eqn:inner_product_half_dens}
    \bracket{\psi_1}{\psi_2} = \int_{\mathcal{D}}\psi^*_1(x)\psi_2(x)\,,
\end{align}
where the extension to higher dimensional Euclidean spaces can be treated similarly.
\begin{remark}\label{rmk:half_den_notation}
    Following \cite{FHT2019}, the standard geometric treatment of complex wave functions is to regard them as half densities, $\psi = \tilde{\psi}\sqrt{d^2x}\in {\rm Den}^{\frac{1}{2}}(\mathcal{D})$. Thus far, we have intentionally surpressed the notation for the basis. In particular, equation \eqref{eqn:inner_product_half_dens} implies that the basis is `hidden' within the notation for $\psi$ and $\psi^*$. The modulus of the complex wave function is defined as
    \begin{equation*}
        |\psi|^2 := \psi \psi^* \in {\rm Den}(\mathcal{D}) \,.
    \end{equation*}
    When we need to suppress the basis entirely, we will use a tilde as in $\tilde{\psi}\sqrt{d^2x}\in {\rm Den}^{\frac{1}{2}}(\mathcal{D})$. In what follows, when vector calculus notation is seen (i.e. $\nabla$, $\p_j$, etc) we mean that this acts on the coefficient of the basis but not the basis itself. When exterior calculus notation is seen (i.e. $\mathcal{L}$ etc), we mean that the operator applied to the whole geometric object including its basis. With these caveats in mind, no confusion should arise when conventional calculus notation also appears.
\end{remark}
% Following \cite{FHT2019}, the standard geometric treatment of complex wave functions are regarded as half densities, i.e. $\psi, \psi^* \in \textrm{Den}^{\frac{1}{2}}(\mathcal{D})$ such that the modulus $|\psi|^2 \in \textrm{Den}(\mathcal{D})$. In basis notation, we have $\psi = \tilde{\psi}\sqrt{d^2 x}$ where $\tilde{\psi}$ is the coefficient of the half-density basis $\sqrt{d^2 x}$. For ease of notation, we shall suppress the basis and work with the notation $\psi$ to denote the product of the coefficients and basis.\\ \\
The linear Schr\"odinger equation in terms of the Hamiltonian operator $\widehat{H}$ is the Euler-Lagrange equation of \eqref{eq:DF var prin},
\begin{align}
    i\hbar\p_t\psi = \widehat{H}\psi\,.
\end{align}
By considering the Hamiltonian functional $H(\psi, \psi^*) := \scp{\psi}{\widehat{H}\psi} =: H[\psi]$, Schr\"odinger's equation can be cast into canonical Hamiltonian form as 
\begin{align}
    i\hbar\p_t\psi = \frac{\delta H}{\delta \psi^*}\,,
\end{align}
where the normalisation for the canonical Poisson brackets is taken as $\{\psi(x), \psi^*(x')\} = -\frac{i}{\hbar}\delta(x-x')$ \footnote{A factor of $\frac{1}{2}$ has been introduced to the canonical Poisson structure of $(\psi, \psi^*)$ relative to reference \cite{FHT2019}.}. Similarly, the NLS equation \eqref{det-NLS} may be derived from the Hamiltonian functional
\begin{align}
    H[\psi,\psi^*] =  \frac12\int_{\mathcal{D}} |\nabla \tilde\psi|^2 + \kappa|\tilde\psi|^4 \, d^2x \,, \label{det-NLS-Ham}
\end{align}
where the notation $\tilde\psi$, as in Remark \ref{rmk:half_den_notation} has been applied to ensure that spatial integration is properly defined.
In 1D, the NLS equation is a completely integrable Hamiltonian system, with an infinity of conserved quantities that all Poisson commute amongst themselves, \cite{A&S1981}. However, in higher dimensions, the NLS equation conserves only the energy $H[\psi,\psi^*]$ and the two cotangent-lift momentum maps which arise from the invariances of the deterministic Hamiltonian $H[\psi,\psi^*]$ in \eqref{det-NLS-Ham} under constant shifts of phase and translations in space. Let $g_t \in \textrm{Diff}(\mathcal{D})$ a time dependent diffeomorphism which act on $\psi$ by pull-back, the Lie derivative $\mathcal{L}_u \psi$ of $\psi$ by $u := \bu\cdot \nabla \in\mathfrak{X}(\mathcal{D})$ can be calculated in terms of basis functions as
\begin{align}
    \mathcal{L}_u \psi := \frac{d}{dt}\bigg|_{t=0}(g_t^*\psi) = \left(\frac{1}{2}(\p_ju_j + u_j\p_j)\psi\right)\,,
\end{align}
where $g_t$ is the flow of $u$. The diamond operation $\psi_2\diamond\psi_1 \in \mathfrak{X}(\mathcal{D})^*$ for $\psi_1,\psi_2\in\textrm{Den}^{\frac{1}{2}}(\mathcal{D})$ can be calculated using the pairing \eqref{eq:complex pairing} to have
\begin{align}
\begin{split}
    \scp{\psi_2}{\mathcal{L}_u\psi_1} &= \Re\int\psi_2^*\left(\frac{1}{2}(\p_ju_j + u_j\p_j)\psi_1\right) \\
   &= \Re\int-\left(\frac{1}{2}\psi_1\nabla\psi_2^* - \frac{1}{2}
   \psi_2^*\nabla\psi_1\right)\cdot \bu =:\scp{-\psi_2\diamond\psi_1}{u}.
\end{split}
\end{align}
The complex wave amplitude can be written in polar form, $\psi:= \sqrt{N}\exp(i\phi)\sqrt{d^2x}$, in terms of its modulus, $N\,d^2x = |\psi|^2$, and phase, $\phi$. The cotangent lift momentum map associated with the action of diffeomorphisms can be easily derived from the application of Noether's theorem \cite{M1927}, and written succinctly in terms of our new variables, $(N,\phi)$, as
\begin{align}
    J := \mb{J}  \cdot d\bx \otimes d^2x = \hbar\Im(\psi^*\nabla\psi)\cdot d\bx = \hbar N\nabla\phi\cdot d\bx\,.
\label{hbar-J}
\end{align}
Here, $N\,d^2x \in \textrm{Den}(\mathcal{D})$ and $\phi \in \mathcal{F}(\mathcal{D})$, which forms the cotangent bundle $T^*\mathcal{F}(\mathcal{D}) \simeq \textrm{Den}(\mathcal{D})\times \mathcal{F}(\mathcal{D})$.
This implies that $J := \mb{J} \cdot d\bx \otimes d^2x \in \Lambda^1(\mathcal{D})\otimes \text{Den}(\mathcal{D})$ is the also the cotangent lift momentum map on $T^*\mathcal{F}(\mathcal{D})$.
Under similar consideration as the case of invariance of translation in space, the invariance of the Hamiltonian to constant phase shift gives the $\varphi \in S^1$ action on $\psi$, given by $\psi \rightarrow e^{i\varphi}\psi$ gives the momentum map $N d^2x = |\psi|^2$. The Hamiltonian functional in \eqref{det-NLS-Ham} can be transformed into 
\begin{align}
H[\phi, N] = \frac{1}{2} \int_{\mathcal{D}} N |\nabla \phi|^2  + |\nabla \sqrt{N}|^2 + \kappa N^2 \, d^2x\,,
\label{Ham-Nphi}
\end{align}
where the Poisson bracket are $\{N, \phi\} = \frac{1}{\hbar}$. The NLS dynamics can be written in $(N, \phi)$ variables as 
\begin{align}
    \begin{split}
    \hbar\p_t \phi &= \big\{\phi, H[\phi, N] \big\} = -\,\frac{\delta H}{\delta N} = -\left(\frac{1}{2}|\nabla \phi|^2 + \frac{1}{8}\frac{|\nabla N|^2}{N^2} - \frac{1}{4}\frac{\Delta N}{N} + \kappa N \right)=: -\,\varpi\,,\\
    \hbar\p_t  N &= \big\{N, H[\phi, N] \big\} = \frac{\delta H}{\delta \phi} =
    -\,{\rm div} \big( N\nabla \phi\big) =: -\,{\rm div} \mb{J}\,,
    \end{split}
    \label{Bernoulli_Q_Law}
\end{align}
where $\varpi$ in equation \eqref{Bernoulli_Q_Law} is the Bernoulli function. According to \eqref{Bernoulli_Q_Law}, the NLS probability density $N$ is advected by the velocity $\mb{J}/N=\nabla \phi$ and the equation for the phase gradient $\nabla \phi$ reduces to the NLS version of Bernoulli's law. 
The Hamiltonian in \eqref{det-NLS-Ham} collectivises through the momentum maps $N$ and $\mb{J}$ into 
\begin{align}
    H[\mb{J},N] = \frac12 \int_{\mathcal{D}} \frac{|\mb{J}|^2}{\hbar^2 N} +  |\nabla \sqrt{N}|^2 + \kappa N^2 \, d^2x\,,
    \label{Ham-NJphi}
\end{align}
such that it is a Hamiltonian functional defined on the semi-direct product Lie algebra $\mathfrak{X}^*(\mathcal{D}) \circledS \textrm{Den}(\mathcal{D})$. The Lie-Poisson structure of $(\mb{J}, N) \in \mathfrak{X}^*(\mathcal{D})\circledS \textrm{Den}(\mathcal{D})$ implies the NLS equation can be expressed in matrix operator Lie-Poisson bracket form as
\begin{align}
    \frac{\p}{\p t}
    \begin{bmatrix}
    J_i \\ N
    \end{bmatrix}
    =
    -
    \begin{bmatrix}
    (\p_kJ_i + J_k\p_i) & N\p_i
    \\ \p_kN & 0
    \end{bmatrix}
    \begin{bmatrix}
    \frac{\delta H[\mb{J}, N]}{\delta J_k} = J_k/(\hbar N) =  \phi_{,k}/\hbar
    \\ 
    \frac{\delta H[\mb{J}, N]}{\delta N} = -\frac{|\mb{J}|^2}{2\hbar^2 N^2} + \frac{1}{8}\frac{|\nabla N|^2}{N^2} - \frac{1}{4}\frac{\Delta N}{N} + \kappa N
    \end{bmatrix}
    .
    \label{NJ-matrixPBs}
\end{align}
Noting that the canonical and the Lie-Poisson Hamiltonian structure of the NLS equation in \eqref{Bernoulli_Q_Law} and \eqref{NJ-matrixPBs} respectively, we can apply both side of the reduction pathway shown in Figure \ref{fig:reduction path} to couple the NLS equation to a fluid flow. 

\begin{remark}
In the following considerations, we shall set $\hbar = 1$ for convenience.
\end{remark}

Let us first consider the coupling of the NLS equation in canonical Hamiltonian form \eqref{Bernoulli_Q_Law} to an inhomogeneous and incompressible Euler's fluid. The fluid Lagrangian $L_{iE}(u, D, \rho)$ is described by the fluid velocity $\bu \cdot \nabla \in \mathfrak{X}(\mathcal{D})$, advected volume density $D d^2x = (D_0 d^2x_0)g^{-1} \in \text{Den}(\mathcal{D})$ and advected, spatially inhomogeneous buoyancy $\rho = \rho_0 g^{-1}\in \mathcal{F}(\mathcal{D})$. The coupling of NLS equation to the inhomogeneous Euler fluid is through the following Hamilton's principle in the form of \eqref{eq:metamorphsis HP var partial Ham},
\begin{equation}
\begin{aligned}
0 &= \delta S = \delta\int_a^b\int_{\mathcal{D}} \ell_{iE}(u,D,\rho) - \langle N, \p_t\phi + \mb{u}\cdot\nabla\phi \rangle - H_{NLS}(N,\phi)
\\
&= \delta\int_a^b\int_{\mathcal{D}} \frac{D\rho}{2} |\mb{u}|^2 - p(D-1) - \mb{u}\cdot  N\nabla\phi - 
 N \p_t\phi - \frac{1}{2} \left( N|\nabla\phi|^2 +  |\nabla \sqrt{N}|^2 + \kappa N^2 \right)
\,d^2x\,dt \,,
\label{HP-NLS-A-Eul-Lag}
\end{aligned}
\end{equation}
where the constrained variations are $\delta u^L = \p_t w - \ad_{u^L} w$, $\delta (D d^3x) = -\mathcal{L}_{w} (Dd^3x) =- {\rm div}(D \mb{w})d^3x$ and $\delta b = -\mathcal{L}_{w} = - \mb{w}\cdot\nabla b$; the arbitrary variations are $\mb{w}\cdot \nabla = w \in \mathfrak{X}(\mathcal{D})$, $\delta N$ and $\delta \phi$ which vanish at at endpoints. In the action \eqref{HP-NLS-A-Eul-Lag}, it is evident that we have coupled the Lagrangian for an inhomogeneous Euler fluid, $\ell_{iE}$, to the Hamiltonian for the NLS equation, $H_{NLS}$, in the manner proposed in equation \eqref{eq:metamorphsis HP var partial Ham}. Namely, we have chosen $G = \Diff(\mathcal{D})$, $V^* = \mathcal{F}(\mathcal{D})\oplus \text{Den}(\mathcal{D})$ and $Q = \mathcal{F}(\mathcal{D})$ in \eqref{eq:metamorphsis HP var partial Ham}. 
Here, the variables $N$ and $\phi$ are pushed-forward by the flow generated by $u$.
The modified canonical Hamiltonian equations for $(N,\phi)$ arising from Hamilton's principle \eqref{HP-NLS-A-Eul-Lag} are
\begin{align}
\begin{split}
\p_t N + {\rm div} \big(N(\mb{u} + \nabla\phi)\big) &=  0
\,,\\
\p_t\phi + \mb{u}\cdot \nabla\phi &=  -\varpi 
\,.
\end{split}
\label{eps-modELeqns-NLS}
\end{align}
Thus, the evolution equations for the Eulerian wave variables $(N,\phi)$ in \eqref{eps-modELeqns-NLS} keep their form as canonical Hamilton's equation forms with the added effects of `Doppler-shifts' by the fluid velocity $\mb{u}$. 
The modified Euler-Poincar\'e equations that arise from Hamilton's principle in \eqref{HP-NLS-A-Eul-Lag} are
\begin{align}
\begin{split}
\big( \p_t + {\cal L}_u \big)\Big( \Big( D\rho\mb{u} - N\nabla\phi\Big) \cdot d\mb{x}\otimes d^2x \Big) 
&= \left(D \nabla \Big(\frac{\rho}{2} |\mb{u}|^2 %+ \rho\mb{u}\cdot \nabla\phi 
- p \Big)
-  \Big(\frac{D}{2} |\mb{u}|^2 %+ D\mb{u}\cdot \nabla\phi  
\Big)d\nabla\rho\right)\cdot d\mb{x}\otimes d^2x
\,,
\end{split}
\label{eps-modELeqnswave}
\end{align}
along with the NLS equations in \eqref{eps-modELeqns-NLS} and the advection equations
\begin{align} 
\begin{split}
\big( \p_t + {\cal L}_u \big)\rho &= \p_t\rho  + \mb{u}\cdot\nabla \rho =  0
\,,\\
\big( \p_t + {\cal L}_u \big)\big( D\, d^2x\big) &= \big(\p_t D + {\rm div} (D\mb{u})\big)\, d^2x = 0
\,,\quad D=1 \Longrightarrow {\rm div} \mb{u} = 0
\,,
\end{split}
\label{eps-advecteqns}
\end{align}
in which preservation of the constraint $D=1$ requires divergence-free flow velocity, ${\rm div} \mb{u} = 0$. Then equations \eqref{eps-modELeqns-NLS} with \eqref{eps-modELeqnswave} imply 
\begin{align}
    \big( \p_t + {\cal L}_u \big)\left(  D\rho \mb{u} \cdot d\mb{x} \otimes d^2x\right) 
    &= \left(D \nabla \Big(\frac{\rho}{2} |\mb{u}|^2 - p  \Big) -  \Big(\frac{D}{2} |\mb{u}|^2 \Big)\nabla \rho 
    - \textrm{div}(N\nabla\phi)\nabla \phi - N \nabla \varpi\right)\cdot d\mb{x}\otimes d^2x \,.
    \label{eps-modELeqnswave-tangled}
\end{align}
The equations \eqref{eps-modELeqnswave-tangled}, \eqref{eps-advecteqns} and \eqref{eps-modELeqns-NLS} are exactly in the general form \eqref{eq:metamorphsis LP eq}. 
The general result in equation \eqref{eq:KC untangled} yields the following Kelvin-Noether theorem for the total Hamilton's principle for NLS waves on a free fluid surface in equation \eqref{HP-NLS-A-Eul-Lag},
\begin{equation}
    \frac{d}{dt}\oint_{c(\bs{u})}\underbrace{ \Big(\bs{u} - \frac{N\nabla \phi}{D\rho}\Big) \cdot\, d\bs{x} }_{\hbox{`Momentum shift'}}
    = \oint_{c(\bs{u})} \left(
    \nabla\left( \frac{|u|^2}{2} \right)- \frac{1}{\rho}\nabla p \right)\cdot d\mb{x}
    \,.
\label{HP1-KN1-NLS}
\end{equation}
Equation \eqref{eps-modELeqnswave-tangled} yields the separated Kelvin-Noether equations as in \eqref{eq:KC tangled},
\begin{align}
\begin{split}    
     \frac{d}{dt}\oint_{c(\bs{u})}   \mb{u} \cdot d\mb{x}  
    &=  \oint_{c(\bs{u})} \left(\nabla \left(\frac{|u|^2}{2}\right) - \frac{1}{\rho}\nabla p\right)\cdot d\mb{x} - \oint_{c(\bs{u})}
    \underbrace{\frac{1}{D\rho}\left(\textrm{div}(N\nabla\phi)\nabla\phi + N\nabla \varpi\right)
    \cdot d\mb{x}}_{\hbox{Non-inertial force}}
\,,
\\ 
\frac{d}{dt}\oint_{c(\bs{u})} \frac{1}{D\rho}\, (N\nabla\phi)\cdot d\mb{x}
&= -\,
\oint_{c(\bs{u})} \,\frac{1}{D\rho}\left(\textrm{div}(N\nabla\phi)\nabla\phi + N\nabla\varpi\right)\cdot d\mb{x}
\,,
\\&= -\,
\oint_{c(\bs{u})} \,\frac{1}{D\rho}\Bigg(\p_j\big(N\phi^{\,,\,j}\phi_{,\,k}\big)dx^k
- \frac{N}{4} \nabla\left( \frac{|\nabla N|^2}{2N^2} - \frac{\Delta N}{N} + 4\kappa N 
\right)\cdot d\mb{x}
\Bigg)\,,
\end{split}
\label{HP1-KN2-NLS}
\end{align}
where $\varpi$ is again the Bernoulli function in equation \eqref{Bernoulli_Q_Law}. The stress tensor $T^j_k:=N\phi^{\,,\,j}\phi_{,\,k}$ in the last equation mimicks the corresponding stress tensor in the evolution of the Berry curvature in quantum hydrodynamics; see equation (106) in \cite{FHT2019}.
\begin{remark}
Upon comparing the unified and separated Kelvin circulation equations in \eqref{HP1-KN1-NLS} and \eqref{HP1-KN2-NLS}, respectively, one sees that: \\
(1)  In \eqref{HP1-KN1-NLS} the standard Kelvin circulation theorem for an inhomogeneous planar Euler flow holds in the absence of waves. Thus, the fluid flow does not create waves.
\\ (2) In \eqref{HP1-KN2-NLS} the first equation of the separated Kelvin theorem shows that the Kelvin circulation theorem for an inhomogeneous planar Euler flow has an additional source in the presence of waves. Thus, one sees that the waves can create circulatory fluid flow.
\end{remark}
%\todo[inline,color=yellow]{DH: Recall our question about ``back-reaction'' or ``backflow'' in the Introduction for this system.}
In terms of the fluid momentum density $\mb{m} := D\rho \mb{u}$ with fluid transport velocity $\mb{u}$, the  Hamiltonian for NLS wave-current system dynamics is written as
\begin{align}
    H_{m}[\mb{m},D,\rho,\phi, N] = \int_{\mathcal{D}} \frac{|\mb{m}|^2}{2D\rho}  + p(D-1)
    + \frac12 \Big(N |\nabla \phi|^2  + |\nabla \sqrt{N}|^2 + \kappa N^2 \Big) \, d^2x\,.
    \label{Ham-mNphi}
\end{align}
The dynamics of the current-coupled NLS system may then be written in Lie-Poisson bracket form as
\begin{align}
    \frac{\p}{\p t}
    \begin{bmatrix}
    m_i \\ D \\ \rho \\ \phi \\ N
    \end{bmatrix}
    = -
    \begin{bmatrix}
    (\p_k m_i + m_k\p_i) & D\p_i&  -\rho_{,i} & -\phi_{,i}& N\p_i
    \\ \p_k D & 0 & 0 & 0 & 0
    \\ \rho_{,k} & 0 & 0 & 0 & 0
    \\ \phi_{,k} & 0 & 0 & 0 & 1
     \\ \p_k N & 0 & 0 & -1 & 0
    \end{bmatrix}
    \begin{bmatrix}
    \frac{\delta H_{m}}{\delta m_k} = u_k 
    \\
    \frac{\delta H_{m}}{\delta D} = -\frac{|\mb{m}|^2}{2 D^2\rho}
    \\
    \frac{\delta H_{m}}{\delta \rho} = -\frac{|\mb{m}|^2}{2 D\rho^2}
    \\
    \frac{\delta H_{m}}{\delta \phi} = -\textrm{div}(N\nabla \phi)
    \\ 
    \frac{\delta H_{m}}{\delta N} = \varpi
    \end{bmatrix}
    ,
    \label{Nphi-eqns}
\end{align}
where the Bernoulli function $\varpi$ is given in equation \eqref{Bernoulli_Q_Law}. By taking the untangling map and writing the Hamiltonian \eqref{Ham-mNphi} in terms of the total momentum $\mb{M}: = \mb{m} - N\nabla \phi$, we have the following Hamiltonian 
\begin{align}
    H_{HP}[\mb{M},D,\rho,\phi, N] = \int_{\mathcal{D}} \frac{|\mb{M} + N\nabla \phi|^2}{2D\rho}  + p(D-1)
    + \frac12 \Big(N |\nabla \phi|^2  + |\nabla \sqrt{N}|^2 + \kappa N^2 \Big) \, d^2x\,,
    \label{Ham-mNphi-tot}
\end{align}
and the untangled Poisson structure
\begin{align}
    \frac{\p}{\p t}
    \begin{bmatrix}
    M_i \\ D \\ \rho \\ \phi \\ N
    \end{bmatrix}
    = -
    \begin{bmatrix}
    (\p_k M_i + M_k\p_i) & D\p_i&  -\rho_{,i} & 0& 0
    \\ \p_k D & 0 & 0 & 0 & 0
    \\ \rho_{,k} & 0 & 0 & 0 & 0
    \\ 0 & 0 & 0 & 0 & 1
     \\ 0 & 0 & 0 & -1 & 0
    \end{bmatrix}
    \begin{bmatrix}
    \frac{\delta H_{HP}}{\delta M_k} = \frac{\delta H_{m}}{\delta m_k} = u_k 
    \\
    \frac{\delta H_{HP}}{\delta D} = -\frac{|\mb{M}+N\nabla \phi|^2}{2 D^2\rho} = \frac{\delta H_{m}}{\delta D}
    \\
    \frac{\delta H_{HP}}{\delta \rho} = -\frac{|\mb{M}+N\nabla \phi|^2}{2 D\rho^2} = \frac{\delta H_{m}}{\delta \rho}
    \\
    \frac{\delta H_{HP}}{\delta \phi} = -\textrm{div}(N(\nabla \phi + \mb{u})) = -\textrm{div}(N\mb{u}) + \frac{\delta H_{m}}{\delta \phi}
    \\ 
    \frac{\delta H_{HP}}{\delta N} = \varpi+\mb{u}\cdot\nabla\phi = \mb{u}\cdot\nabla\phi + \frac{\delta H_{m}}{\delta N}
    \end{bmatrix}\,.
    \label{Nphi-untangled-eqns}
\end{align}

The transformation to the Lie-Poisson wave variables $(\mb{J}, N)$, the canonical Hamiltonian \eqref{Ham-mNphi} transforms to
\begin{align}
    H_{J}[\mb{m},D,\rho,\mb{J},N] =  \int_{\mathcal{D}} \frac{|\mb{m}|^2}{2D\rho}  + p(D-1)
    + \frac{|\mb{J}|^2}{2 N} + \frac12\Big( |\nabla \sqrt{N}|^2 + \kappa N^2\Big) \, d^2x\,,
    \label{Ham-mNJ}
\end{align}
and the corresponding equations in Lie-Poisson bracket form are given by 
\begin{align}
        \frac{\p}{\p t}
    \begin{bmatrix}
    m_i \\ D \\ \rho \\ J_i \\ N
    \end{bmatrix}
    = -
    \begin{bmatrix}
    (\p_k m_i + m_k\p_i) & D\p_i&  -\rho_{,i} & (\p_k J_i + J_k\p_i) & N\p_i
    \\ \p_k D & 0 & 0 & 0 & 0
    \\ \rho_{,k} & 0 & 0 & 0 & 0
    \\ (\p_k J_i + J_k\p_i)  & 0 & 0 & (\p_k J_i + J_k\p_i)  & N\p_i
     \\ \p_k N & 0 & 0 & \p_k N & 0
    \end{bmatrix}
    \begin{bmatrix}
    \frac{\delta H_{J}}{\delta m_k} = u_k 
    \\
    \frac{\delta H_{J}}{\delta D} = -\frac{|m|^2}{2 D^2\rho}
    \\
    \frac{\delta H_{J}}{\delta \rho} = -\frac{|m|^2}{2 D\rho^2}
    \\
    \frac{\delta H_{J}}{\delta J_k} = J_k/N
    \\ 
    \frac{\delta H_{J}}{\delta N} = \Tilde{\varpi}
    \end{bmatrix}
    \, \label{eq:mDrhoNJ-matrixPBs}
\end{align}
where $\Tilde{\varpi} := -\frac{|\mb{J}|^2}{2N^2} + \frac{1}{8}\frac{|\nabla N|^2}{N^2} - \frac{1}{4}\frac{\Delta N}{N} + \kappa N$. 
In transforming the wave variables from $(\phi,N)$ to $(\mb{J},N)$ the canonical two-cocyle between $(\phi, N)$ has been transformed into a generalised cocycle in $(\mb{J}, N)$ which takes the form of the Lie-Poisson bracket on the dual of $\mathfrak{X}(\mathcal{D})\circledS \mathcal{F}(\mathcal{D})$ and it is a particular case of the general form of Lie-Poisson system with generalised co-cycle \eqref{G Q metamorphsis PB}.
Specifically, the Poisson bracket \eqref{eq:mDrhoNJ-matrixPBs} is a standard Lie-Poisson bracket on the dual of the Lie algebra 
\begin{align}
\mathfrak{X}_1(\mathcal{D})\circledS\big((\mathfrak{X}_2(\mathcal{D}) \circledS \mathcal{F}(\mathcal{D}))\oplus \mathcal{F}(\mathcal{D})\oplus \textrm{Den}(\mathcal{D})\big)
\,,\end{align}
where the corresponding semidirect-product Lie group is 
\begin{align}
\textrm{Diff}_1(\mathcal{D})\circledS\big((\textrm{Diff}_2(\mathcal{D})\circledS\mathcal{F}(\mathcal{D}))\oplus\mathcal{F}(\mathcal{D})\oplus\textrm{Den}(\mathcal{D})\big)
\,.
\end{align}
Equation \eqref{eq:mDrhoNJ-matrixPBs} yields a modified version of separated Kelvin-Noether theorem, namely,
\begin{align}
    \begin{split}
    \frac{d}{dt}\oint_{c(\bs{u})} \mb{u} \cdot d\mb{x}  
    &=  \oint_{c(\bs{u})} \left(\nabla \left(\frac{|u|^2}{2}\right) - \frac{1}{\rho}\nabla p\right)\cdot d\mb{x} \\
    &\qquad - \oint_{c(\bs{u})}
    \underbrace{\frac{1}{D\rho}\left( \frac{\mb{J}}{N} \cdot \nabla \mb{J} + J_k\nabla \left(\frac{J_k}{N}\right) + \mb{J}\textrm{div}(\mb{J}/N) + N\nabla \Tilde{\varpi}\right)\cdot d\mb{x}
    }_{\hbox{Non-inertial force}}
    \,,
    \\ 
    \frac{d}{dt}\oint_{c(\bs{u})} \frac{1}{D\rho}\, \mb{J}\cdot d\mb{x}
    &= -\, \oint_{c(\bs{u})} \,\frac{1}{D\rho}\left( \frac{\mb{J}}{N} \cdot \nabla \mb{J} + J_k\nabla \left(\frac{J_k}{N}\right) + \mb{J}\textrm{div}(\mb{J}/N) + N\nabla \Tilde{\varpi}\right)\cdot d\mb{x}
    \,,
    \end{split}
\end{align}

\begin{remark}[Coupling to complex half densities]
For completeness, let us consider Hamilton's principle for coupling the inhomogenous Euler's equation to the NLS equations in the complex wave function variables $(\psi, \psi^*)$ \eqref{det-NLS}, it reads,
\begin{align}
0 = \delta S &= \delta\int_a^b \int_\mathcal{D}\left(\frac{D\rho}{2} |\mb{u}|^2 
- p(D-1) - \mb{u}\cdot \Im(\psi^*\nabla\psi)\right)\,d^2x + \scp{\psi}{i\p_t\psi} - H[\psi, \psi^*] \,dt \,,
\label{HP-NLS-A-Eul-Lag-psi}
\end{align}
where $H[\psi,\psi^*]$ is the NLS Hamiltonian in terms of $(\psi, \psi^*)$, defined in \eqref{det-NLS-Ham}.
The canonical equations for complex wave function $\psi$ can then be calculated to be 
\begin{align}
    i\hbar\left(\p_t + \mathcal{L}_u\right)\psi := i\hbar\left(\p_t + \frac{1}{2}(\p_ju^j + u^j\p_j)\right)\psi = -\frac{1}{2}\triangle\psi + \kappa|\psi|^2\psi\,.
    \label{half-dens}
\end{align}
Just as the current boosts the scalar phase $\phi$ and density $N d^2x$ by the Lie derivative in equation \eqref{eps-modELeqnswave}, the half density $\psi\sqrt{d^2x}$ is also boosted by the Lie derivative with respect to the current velocity vector field $u$ in equation \eqref{half-dens}.
\end{remark}
\begin{remark}[Coupling NLS to mesoscale QG motion]
Coupling of NLS to homogeneous $(\rho=1)$ mesoscale QG motion can be accomplished by modifying the reduced Lagrangian in \eqref{HP-NLS-A-Eul-Lag-psi} to include rotation and quasigeostrophic balance, as follows \cite{HZ1998,Z2018}
\begin{align}
0 = \delta S &= \delta\int_a^b \int_\mathcal{D} \bigg[
\frac{D}{2}\Big( \mb{u}\cdot \big(1-\mc{F}\Delta^{-1}\big)\mb{u} + \mb{u}\cdot \mb{R}(\mb{x})\Big)
- p(D-1) 
\\&\hspace{2cm}- \mb{u}\cdot \Im(\psi^*\nabla\psi)\bigg]\,d^2x + \scp{\psi}{i\p_t\psi} - H[\psi, \psi^*] \,dt 
\,.
\label{HP-NLS-A-QG-Lag-psi}
\end{align}
Here, $\mc{F}$ is the rotational Froude number and $\mb{R}(\mb{x})$ is the prescribed vector potential for the Coriolis parameter. The derivation of the equations of motion and Hamiltonian formulation can be accomplished by combining the calculations above with those in \cite{HZ1998,Z2018} to accommodate rotation and quasigeostrophy.

\end{remark}

\section{Numerical simulations}\label{sec:4-Numerics}
\begin{figure}[h!]
    \centering
    \begin{subfigure}[b]{0.40\textwidth}
		\centering
		\includegraphics[width=\textwidth]{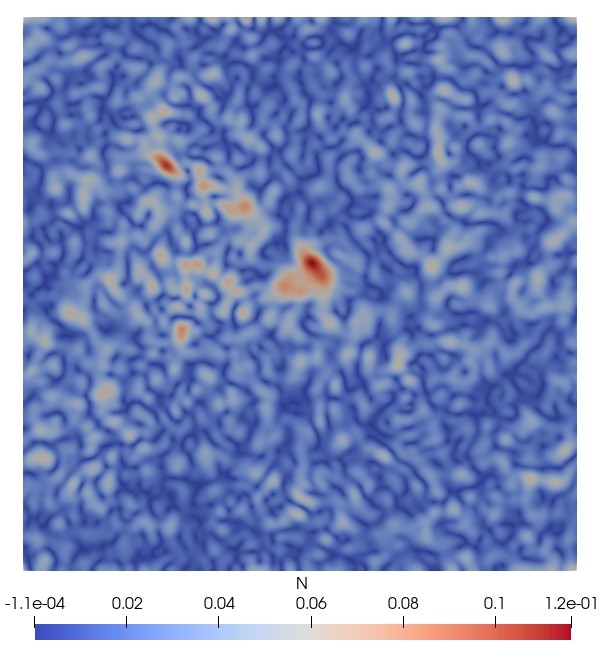}
	\end{subfigure}
	\begin{subfigure}[b]{0.40\textwidth}
		\centering
		\includegraphics[width=\textwidth]{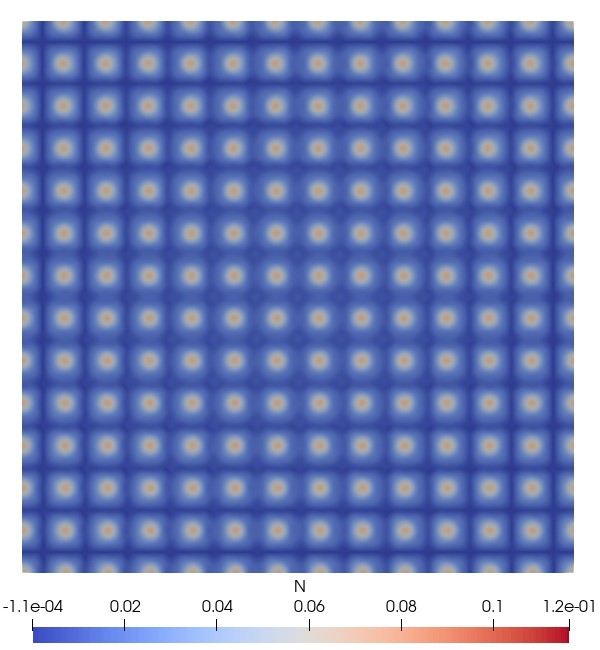}
	\end{subfigure}
	\caption{These are the $512^2$ snapshot of the wave amplitudes $N:= a^2+b^2$ from the numerical simulating the Euler coupled NLS equation \eqref{eq:PV-NLS-eqns}(left) and numerical simulations of the uncoupled NLS equations (right) at time $t=30$. The initial conditions for $Q_F$ and $\rho$ are obtained following a spin-up period of the inhomogeneous Euler equations without waves. As seen in the right hand panel, the uncoupled NLS equation produced a `Gingham' pattern due to the boundary conditions and the spatial symmetry of the initial conditions. However, when coupled to the `mixing' flow of the inhomogeneous Euler's equation, the spatial coherence of $N$ is distorted as seen in the left hand panel. However, it still retains the localisation of the patterns as local regions of high densities usually have filaments of zero densities as boundaries.}
    \label{fig:current on waves N}
\end{figure}

In preparation for the numerical simulations of the coupled non-homogeneous Euler coupled NLS equations \eqref{eq:mDrhoNJ-matrixPBs} in 2D, as discussed in Section \ref{subsec:NLS}, let us consider the equation in terms of the real and imaginary parts of $\psi$, namely $a$ and $b$ such that $\psi := a+ib$. This particular change of variables is done for ease of implementation of the numerical solver. Inserting these relations into the action \eqref{HP-NLS-A-Eul-Lag-psi} gives 
\begin{align}
\begin{split}
    0 = \delta S &= \delta\int_a^b \int_\mathcal{D}\frac{D\rho}{2} |\mb{u}|^2 
    - p(D-1) + \hbar\left(b\left(\p_t + \mb{u}\cdot \nabla\right)a - a\left(\p_t + \mb{u}\cdot \nabla\right)b\right)\\
    &\qquad \qquad - \frac{1}{2}\left(|\nabla a|^2 + |\nabla b|^2 + \kappa\left(a^2 + b^2\right)^2\right) \,d^2x\,dt \,,
\end{split} \label{HP-NLS-A-Eul-Lag-ab}
\end{align}
The NLS momentum map in terms of $a, b$ can be computed as $\mb{J}(a,b) := \hbar(a\nabla b - b\nabla a)$ and we have the equation to solve as 
\begin{align}
    \begin{split}
        & \left(\p_t + \mathcal{L}_{\mb{u}}\right)\left(\left(D\rho\mb{u} - \mb{J}(a,b)\right)\cdot d\mb{x}\right) = Dd\left(\frac{\rho}{2}|\mb{u}|^2 - p\right) - \frac{D}{2}|\mb{u}|^2d\rho\,, \\
        & \p_t\rho + \mb{u}\cdot\nabla\rho = 0\,,\quad \p_t D + \textrm{div}(D\mb{u}) = 0 \,,\quad D = 1 \Rightarrow \textrm{div}(\mb{u}) = 0\,,\\
        & \p_t a + \mathcal{L}_{\mb{u}} a = -\frac{1}{2}\Delta b + \kappa\left(a^2+b^2\right)b\,,\\
        &\p_t b + \mathcal{L}_{\mb{u}} b = \frac{1}{2}\Delta a - \kappa\left(a^2+b^2\right)a\,,
    \end{split}
\end{align}
where we have again set $\hbar = 0$ for convience. In 2D, one can cast the equation into stream function and vorticity form by defining fluid and wave potential vorticities as follows
\begin{align}
    Q_F\,d^2x := d\left(\rho \mb{u}\cdot d\mb{x}\right) = \textrm{div}(\rho\nabla \Psi)\,, \quad Q_W\,d^2x := d\left(\mb{J}(a,b)\cdot d\mb{x}\right) = 2\hbar J(a,b)d^2x\,,
\end{align}
where $\Psi$ is the stream function, $\mb{u} = \nabla^\perp \Psi$ and the Jacobian operator $\mathcal{J}$ is defined by $\mathcal{J}(f,h) = \p_x f \p_y h - \p_y f \p_x h$ for arbitrary smooth functions $f,h$. In these variables, the Euler-NLS equations take the following form,
\begin{align}
\begin{split}
    &\p_t (Q_F - Q_W) + \mathcal{J}(\Psi, Q_F - Q_W) = \frac{1}{2}\mathcal{J}(\rho, |\mb{u}|^2)\,,\\
    &\p_t Q_W + \mathcal{J}(\Psi, Q_W) = 2\mathcal{J}\left(-\frac{1}{2}\Delta b + \kappa(a^2 + b^2)b, b\right) + 2\mathcal{J}\left(a, \frac{1}{2}\Delta a - \kappa(a^2 + b^2)a\right)\,,\\
    &\p_t \rho + \mathcal{J}(\Psi, \rho) = 0\,,\\
    &\p_t a + \mathcal{J}(\Psi, a) = -\frac{1}{2}\Delta b + \kappa\left(a^2+b^2\right)b\,,\\
    &\p_t b + \mathcal{J}(\Psi, b) = \frac{1}{2}\Delta a - \kappa\left(a^2+b^2\right)a\,.
\end{split} \label{eq:PV-NLS-eqns}
\end{align}
Our implementation of the inhomogeneous Euler coupled NLS equations \eqref{eq:PV-NLS-eqns} used the finite element method (FEM) for the spatial variables. The FEM algorithm we used is implemented using the Firedrake \footnote{\url{https://firedrakeproject.org/index.html}} software. In particular, for \eqref{eq:PV-NLS-eqns} we approximated the fluid potential vorticity $Q_F$, buoyancy $\rho$ using a first order discontinuous Galerkin finite element space. As these quantities are discontinuous on cell boundaries, we have also introduced upwinding in the direction of $\bu$ for the equation of $Q_F$ and $\rho$ for stability of the numerical method.
The real and imaginary parts of the complex wave function, $a$ and $b$, and the stream function $\Psi$ are approximated using a first order continuous Galerkin finite element space. For the time integration, we used the third order strong stability preserving Runge Kutta method \cite{Got05}.
In the numerical examples, we demonstrate numerically the effects of currents on waves and the effects of waves on currents by considering two runs of the 2D inhomogeneous Euler coupled NLS equations \eqref{eq:PV-NLS-eqns} with the following parameters. The domain is $[0,50]^2$ at a resolution of $512^2$. The boundary conditions are periodic in the $x$ direction, homogeneous Dirichlet for $\Psi$, homogeneous Neumann for $a$ and $b$ in the $y$ direction.
To see the effects of the waves on the currents, the procedure was divided into two stages. The first stage was performed on the inhomogenous Euler's equations for $T_{spin} = 100$ time units starting from the following initial conditions
\begin{align}
\begin{split}
    Q_F(x,y, 0) & = \sin(0.16\pi x)\sin(0.16\pi y) + 0.4\cos(0.12\pi x)\cos(0.12\pi y) + 0.3\cos(0.2\pi x)\cos(0.08\pi y) +\\ & \qquad 0.02\sin(0.04\pi y) + 0.02\sin(0.04\pi x)\,,\\
    \rho(x,y, 0)& = 1 + 0.2\sin(0.04\pi x)\sin(0.04\pi y)\,.
\end{split}
\end{align}
The purpose of the first stage was to allow the fluid system to spin up to a statistically steady state without influences from the wave dynamics. The PV and buoyancy variables at the end of the initial spin-up period are denoted as $Q_{spin}(x,y) = Q_F(x,y,T_{spin})$ and $\rho_{spin}(x,y) = \rho(x,y,T_{spin})$. In the second stage, the full simulations including the wave variables were run with the initial conditions for the fluid variables being the state achieved at the end of the first stage. To start the second stage for \eqref{eq:PV-NLS-eqns}, wave variables were introduced with the following initial conditions
\begin{align}
    \begin{split}
        a(x,y,0) & = \exp(-((x-25)^2+(y-25)^2))\,,\quad b(x,y,0) = 0\,, \quad \kappa = \frac{1}{2}\,,\\
        Q_F(x,y,0) &= Q_{spin}(x,y)\,,\quad \rho(x,y,0) = \rho_{spin}(x,y)\,.\label{eq:PV stage2 init cond}
    \end{split}
\end{align}
For comparison, we also consider the numerical simulations of the 2D NLS equation without coupling to the inhomogenous Euler equation. The uncoupled NLS equations in the $a$ and $b$ variables are simply the last two equations of \eqref{eq:PV-NLS-eqns} with $\Psi = 0$. From the same initial condition \eqref{eq:PV stage2 init cond}, the snapshots at $t = 30$ of the coupled and uncoupled equations are shown in Figure \ref{fig:current on waves N}.
% \begin{figure}[H]
%     \centering
%     \begin{subfigure}[b]{0.40\textwidth}
% 		\centering
% 		\includegraphics[width=\textwidth]{NLS_defocus_Ninit.jpeg}
% 	\end{subfigure}
% 	\begin{subfigure}[b]{0.40\textwidth}
% 		\centering
% 		\includegraphics[width=\textwidth]{NLS_defocus_Qinit.jpeg}
% 	\end{subfigure}
%  \caption{}
%     \label{fig:current on waves N init}
% \end{figure}

\begin{figure}[h!]
    \centering
    \begin{subfigure}[b]{0.40\textwidth}
		\centering
		\includegraphics[width=\textwidth]{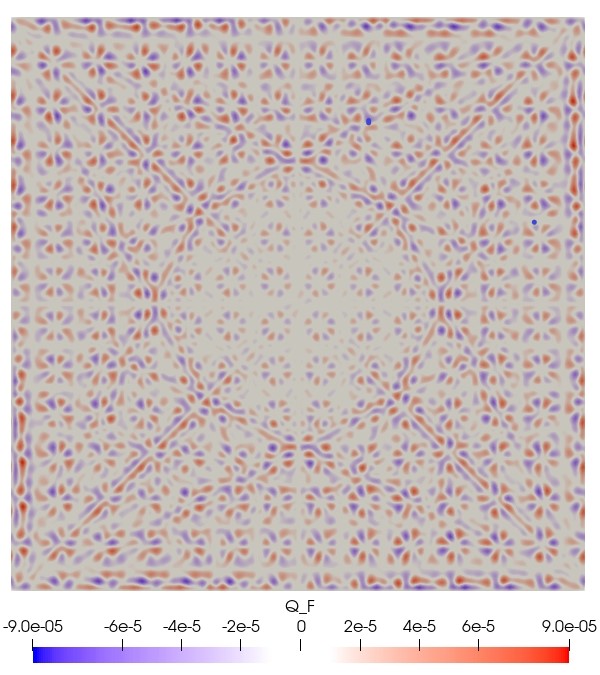}
	\end{subfigure}
	\begin{subfigure}[b]{0.40\textwidth}
		\centering
		\includegraphics[width=\textwidth]{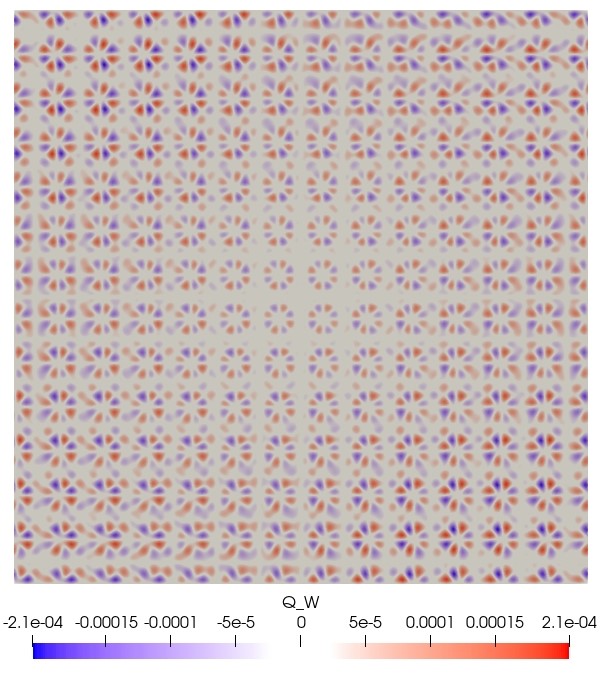}
	\end{subfigure}
	\caption{These are the $512^2$ snapshot of the fluid PV $Q_F$ (left) and wave PV $Q_W$ (right) snapshots at time $t=30$ of the numerical simulation of \eqref{eq:PV-NLS-eqns} with the zero fluid PV initial conditions \eqref{eq:wave on current init cond}. From the right hand panel, one sees the $Q_W$ field form a coherent spatial pattern similar to wave amplitude $N$ of the uncoupled NLS simulation in the left panel of Figure \ref{fig:current on waves N}. The left hand panel is the $Q_F$ generated by $Q_W$. The overall patterns of $Q_F$ is reminiscent of $Q_W$, however, $Q_F$ also shows signs of `mixing' by the fluid since the generated fluid PV will interact with buoyancy to generate circulation. Note that the magnitude of the $Q_F$ is much smaller than $Q_W$, thus the isolated NLS dynamics is dominant over the advection dynamics which implies the minimal `mixing' in the $Q_W$ field. }
	\label{fig:waves on current PV}
\end{figure}

To show the effect of the waves on the currents, we consider a numerical simulation with the following initial conditions,
\begin{align}
    \begin{split}
        a(x,y,0) & = \exp(-((x-25)^2+(y-25)^2))\,,\quad b(x,y,0) = 0\,, \quad \kappa = \frac{1}{2}\,,\\
        Q_F(x,y,0) &= 0\,,\quad \rho(x,y,0) = \rho_{spin}(x,y)\,.\label{eq:wave on current init cond}
    \end{split}
\end{align}
In \eqref{eq:wave on current init cond}, we have used the same initial condition as in \eqref{eq:PV stage2 init cond} except from the PV $Q_F$ which has been set to zero. With this configuration, any PV excitation generated by the waves can interact with a ``well mixed'' buoyancy field to generate further circulation. Snapshots of the $Q_F$ and $Q_W$ fields are shown in Figure \ref{fig:waves on current PV} for the numerical simulations started from the initial conditions \eqref{eq:wave on current init cond}. 
From Figure \ref{fig:waves on current PV}, we see that the spatial features of $Q_W$ are localised and periodic in both directions with varying densities. $Q_F$ possess similar spatial features as $Q_W$, however, these features are deformed. The deformations are precisely caused by the transport of the generated fluid flow and interaction with the buoyancy field. 

\section{Conclusion and outlook}\label{sec:5-Conclusion}

{\bf Summary.} After reviewing the framework in geometric mechanics for deriving hybrid fluid models in the introduction, section \ref{sec:2-Lagrangian reduction} showed a path for their derivation, section \ref{sec:3-NLS} discusssed examples of the wave mean flow hybrid equations and section \ref{sec:4-Numerics} showed simulations of the hybrid Euler-NLS equations. The hybrid Euler-NLS equations describe boosted dynamics of small-scale NLS subsystems into the moving frame of the large-scale 2D Euler fluid dynamics. The Kelvin-Noether theorem in section \ref{sec:2-Lagrangian reduction} showed that the small-scale dynamics can feed back to create circulation in the large-scale dynamics. Over a short time, this creation of large-scale circulation may be only a small effect, as shown in numerical simulations displayed in Figures 2 and 3 of section \ref{sec:4-Numerics}. Over a long enough time period, though, the small-scale effects may produce a more pronounced effect on the larger scales, especially if the small-scale momentum is continuously driven externally. 

{\bf Waves versus patterns.} NLS is a pattern-forming equation that is associated with several different applications in several different fields, including nonlinear fibre optics dynamics of telecommunication as well as studies of deep water waves. When linear driving and dissipation are introduced, NLS becomes the Complex Ginzburg Landau (CGL) equation, which is another well-known pattern-forming equation, \cite{AK2002,M2000, MS2004}. This class of equations is extremely useful for its universal quality as normal form equations for bifurcations, the onset of instability due to symmetry breaking, and the saturation of instability due to nonlinear effects \cite{MS2004}. 
The utility of CGL universality suggests, in particular, that a dissipative and driven version of the hybrid Euler-NLS equations -- that is, the hybrid Euler Complex Ginzburg--Landau (ECGL) equations -- could be proposed as an elementary model to describe some aspects of air-sea coupling that can be encompassed with only a few parameters.  Computational simulations of this proposition are to be discussed elsewhere in future work.

\subsection{Acknowledgements}
This paper was written in appreciation of the late Hermann Flaschka’s elegant, thoughtful and sometimes humorous contributions to nonlinear mathematics during his marvellous career. We hope that the paper has presented ``do-able examples that reveal something new.'' (Namely, that waves are not always carried passively by the current. Waves  can feed back in the Kelvin theorem to produce circulation of the mean fluid velocity that carries them.)  We are grateful to our friends, colleagues and collaborators for their advice and encouragement in the matters treated in this paper. 
DH especially thanks C. Cotter, C. Tronci, F. Gay-Balmaz, T. S. Ratiu and the late J. E. Marsden for many insightful discussions of corresponding results similar to the ones derived here for WMFI, and in earlier work together in deriving hybrid models of complex fluids, turbulence, plasma dynamics, vertical slice models and the quantum--classical hydrodynamic description of molecules. 
OS was partially supported during the present work by European Research Council (ERC) Synergy grant STUOD -- DLV-856408. DH was partially supported during the present work by European Research Council (ERC) Synergy grant STUOD -- DLV-856408 and Office of Naval Research (ONR) grant (Grant No. N00014-22-1-2082). RH was supported during the present work by an EPSRC scholarship (Grant No. EP/R513052/1). %(And later, ONR grant \dots!)

% \subsection{Rigid body version}

% The equations for $(TSO(3) \times TS^1)/SO(3)$ differ from those for $TSO(3)/SO(3) \times TS^1$.  

% So, should we not compare them for the rigid body? Don’t forget, it’s reduction by left action.

% The Lagrangian for for $TSO(3)/SO(3) \times TS^1$ for a flywheel attached to the 2-axis of the rigid body is treated on page 354 of my Blue Book. Should take a look at Krishnaprasad and Marsden's paper on rigid body (satellite) connected to a flaxible cable. 

% \todo[inline]{RH: Krishnaprasad, P.S., Marsden, J.E. Hamiltonian structures and stability for rigid bodies with flexible attachments. Arch. Rational Mech. Anal. 98, 71–93 (1987). \url{https://doi.org/10.1007/BF00279963}. They had the same Poisson structure as ours in section \ref{sec:Lagrangian reduction} in the tangled form because they made the observation that ``Since the body and attachment can be simultaneously rotated''. They also made the observation that one can remove coupling by making a change of momentum as the untangling map. \\
% I think for the Lagrange gryostat, the action of $SO(3)$ on $S^1$ is apparent in the equation for $\alpha$, it is $\dot{\alpha} = \frac{\ell_2}{J_2} - \Omega_2$ takes the form of a Lie-derivative and the variational derivative of a Hamiltonian $h = \frac{1}{2}\frac{\ell_2^2}{J_2}$. In fact, the Lagrange gyrostat equations are Lagrange Poincar\'e equations, see section $3.13$ of \url{https://doi.org/10.1007/978-1-4939-3017-3}.}

%%%%%%%%%%%%%%%%%%%%%%%%%%%%%%%%%%%%%%%%%%%%%%%%%%%%%%%%%%%%%

\begin{appendices}
\section{Stochastic Hamiltonian wave-current dynamics}\label{appendix: stochastic WCI}

The inclusion of a stochastic representation of uncertainty within the model equations has emerged as an effective method of parametrising information lost during the process of discretising and numerically integrating a partial differential equation governing fluid motion. Through the variational structure, we can include a stochastic noise into the model whilst preserving certain desirable properties \cite{H2015}.

\begin{remark}[Two stochastic systems.]\label{rmk:TwoNoises}
    Since our Lagrangian has dependence on two configuration spaces, $G$ and $Q$, there are two dynamical systems which may be made stochastic. In order to state the equations of motion in their most general form, we will here make both systems stochastic. The fluid system, corresponding to the group $G$, will feature stochastic transport noise in the sense of Holm \cite{H2015}, and the noise in the wave system will be expressed as a stochastic version of Hamilton's canonical equations, pioneered by Bismut \cite{B1981}.
\end{remark}

Following Remark \ref{rmk:TwoNoises}, we incorporate noise into the fluid transport velocity, writing the \emph{stochastic transport} vector field in a compact form as
\begin{equation}\label{eqn:SALT_VF}
    \rmd x_t = u\,dt + \sum_i \xi_i\circ dW_t^i \,,
\end{equation}
where $\{W_t^i\}_{i\in\mathbb{N}}$ are independent and identically distributed (i.i.d.) Brownian processes. The coefficients of the noise terms in the transport velocity corresponding to horizontal currents are expressed as vector fields $\xi_i \in \mathfrak{g}$, whereas for the wave dynamics we will define them as variational derivatives of a family of Hamiltonians, $\{ \bar{h}_i(\pi,n) \}$, which are defined by the choice of the stochastic wave equation. The stochastic time integration in the wave dynamics is performed with respect to a distinct collection of Brownian paths. As such, in the stochastic case, the Hamiltonian can be expressed in terms of the Lagrangian as
\begin{equation}\label{eqn:StochasticHamiltonian}
\begin{aligned}
    h(m,n,\pi,a)\,dt + \sum_i h_i(m,a)\circ dW_t^i + \sum_i \bar{h}_i(n,\pi)\circ d\overline W_t^i &= \langle m,\rmd x_t \rangle + \langle \pi,\nu \rangle \\
    &\quad- \ell(u,n,\nu,a) + \sum_i \bar{h}_i(n,\pi)\circ d\overline W_t^i \,,
\end{aligned}
\end{equation}
where $\{W_t^i,\overline W_t^i \}_{i\in\mathbb{N}}$ are i.i.d. Brownian processes. Thus, the inclusion of the noise terms in the equations of motion can be achieved through exploiting the Poisson structure given in \eqref{eq:metamorphosis Poisson matrix}. The stochastic equations are
\begin{align}
    \rmd\begin{pmatrix}m\\a\\ \pi \\ n \end{pmatrix} =
    -\begin{pmatrix}\ad^*_{\fbox{}}m & \fbox{}\diamond a & \fbox{}\diamond \pi & \fbox{}\diamond n \\
    \mathcal{L}_{\fbox{}}a & 0 & 0 & 0\\
    \mathcal{L}_{\fbox{}}\pi & 0 & 0 & 1\\
    \mathcal{L}_{\fbox{}}n & 0 & -1 & 0\end{pmatrix}
    \begin{pmatrix} \delta h / \delta m \,dt + \sum_i\delta h_i / \delta m \circ dW_t^i \\ {\delta h} / {\delta a}\,dt + \sum_i {\delta h_i}/{\delta a}\circ dW_t^i \\ 
    {\delta h}/{\delta \pi}\,dt + \sum_i {\delta \bar{h}_i}/{\delta \pi}\circ d\overline W_t^i \\ {\delta h}/{\delta n}\,dt + \sum_i {\delta \bar{h}_i}/{\delta n}\circ d\overline W_t^i\end{pmatrix}\,.\label{eq:metamorphosis Poisson matrix stochastic}
\end{align}

By writing the momentum equation in a manner consistent with the untangled Poisson structure \eqref{eq:untangled Poisson matrix}, and the wave equations as they appear in equation \eqref{eq:metamorphosis Poisson matrix stochastic}, we may write the equations of motion as
\begin{equation}\label{eq:Poisson stochastic}
\begin{aligned}
    ( \rmd + {\rm ad}^*_{\rmd x_t} )(m + \pi \diamond n) &= -\frac{\delta h}{\delta a}\diamond a \,dt 
     -\sum_i \frac{\delta h_i}{\delta a}\circ dW_t^i\,,\\
    (\rmd + \mathcal{L}_{\rmd x_t})\pi &= -\frac{\delta h}{\delta n}\,dt - \sum_i\frac{\delta \bar{h}_i}{\delta n}\circ d\overline W_t^i
    \,,\\
    (\rmd + \mathcal{L}_{\rmd x_t}) n &= \frac{\delta h}{\delta \pi}\,dt + \sum_i\frac{\delta \bar{h}_i}{\delta \pi}\circ d\overline W_t^i
    \,,\\
    (\rmd + \mathcal{L}_{\rmd x_t})a &= 0
    \,.
\end{aligned}
\end{equation}
In equations for $\pi$ and $n$, we see on the left hand side that the transport by the currents has been stochastically perturbed (as in equation \eqref{eqn:SALT_VF}). On the right hand side of these equations, we see that the wave motion now has a stochastic Hamiltonian structure. For an example of such a stochastically perturbed Hamiltonian system for water wave dynamics, see \cite{S2022b}.

The first equation in \eqref{eq:Poisson stochastic} implies the following stochastic version of the Kelvin-Noether circulation theorem
\begin{equation}
    \rmd \oint_{c(\rmd x_t)} m + \pi\diamond n = \oint_{c(\rmd x_t)} \frac{\delta\ell}{\delta a}\diamond a \,dt \,.
\end{equation}
The reader should note that this relationship exists by design of the stochasticity \cite{H2015}, rather than by coincidence.

% \begin{align}
% \begin{split}
%     \left(\p_t + \ad^*_u\right)m &=  - \frac{\delta h}{\delta n}\diamond n - \frac{\delta h}{\delta \nu}\diamond \nu - \frac{\delta h}{\delta a}\diamond a\,, \\
%     \left(\p_t + \mathcal{L}_u \right)\pi &= \frac{\delta h}{\delta n}\,,\\
%     \left(\p_t + \mathcal{L}_u \right)n &= \frac{\delta h}{\delta \pi}\,,\\
%     \left(\p_t + \mathcal{L}_u \right)a &= 0\,, \quad\hbox{where}\quad u:=\frac{\delta h}{\delta m}\,,
% \end{split} \label{eq:metamorphsis LP eq}
% \end{align}
% which are the Lie-Poisson equations with cocycles. One can arrange \eqref{eq:metamorphsis LP eq} into Poisson bracket form as 
% \begin{align}
%     \p_t\begin{pmatrix}m\\a\\ \pi \\ n \end{pmatrix} =
%     -\begin{pmatrix}\ad^*_{\fbox{}}m & \fbox{}\diamond a & \fbox{}\diamond \pi & \fbox{}\diamond n \\
%     \mathcal{L}_{\fbox{}}a & 0 & 0 & 0\\
%     \mathcal{L}_{\fbox{}}\pi & 0 & 0 & 1\\
%     \mathcal{L}_{\fbox{}}n & 0 & -1 & 0\end{pmatrix}
%     \begin{pmatrix}\frac{\delta h}{\delta m} = u\\ \frac{\delta h}{\delta a} = -\frac{\delta \ell}{\delta a} \\ \frac{\delta h}{\delta \pi} = \nu \\ \frac{\delta h}{\delta n} = -\frac{\delta \ell}{\delta n}\end{pmatrix}\,.\label{eq:metamorphosis Poisson matrix}
% \end{align}
%
%
%
%

\section{Coupling of Harmonic Oscillations}\label{appendix:SHO}
The variational principle given in equation \eqref{eq:metamorphsis HP var partial Ham} may produce wave activity connected to the fluid motion through the kinematic boundary condition, depending on the choice of wave Hamiltonian and fluid Lagrangian. For example, consider a situation where the two dimensional mean fluid flow is governed by the incompressible Euler equation and the wave-like disturbances around the mean flow are taken to be a field of  harmonic oscillators. In which case, the Lagrangian for the waves is
\begin{equation}
	\ell_w = \frac12\int \dot\zeta^2 - \alpha \zeta^2 \,d^2x \,, 
\end{equation}
for some $\alpha\in\mathbb{R}$. A Legendre transform allows us to determine the Hamiltonian, where the conjugate momentum is
\begin{equation*}
	w = \frac{\delta\ell_w}{\delta \dot\zeta} = \dot\zeta \,,
\end{equation*}
which implies that the wave Hamiltonian is
\begin{equation}
	h_w = \left\langle w , \dot\zeta \right\rangle - \frac12\int \dot\zeta^2 - \alpha \zeta^2 \,d^2x = \frac12\int  w^2 + \alpha\zeta^2 \,d^2x \,.
\end{equation}
Coupling this to the Euler equations, using the approach outlined in Section \ref{sec:2-Lagrangian reduction}, implies an action 
\begin{equation}
\begin{aligned}
	S &= \int\left[\int \frac{D\rho}{2}|\bs{u}|^2 - p(D-1) \,d^2x + \langle w , (\p_t + \mathcal{L}_u)\zeta \rangle - \frac12\int w^2 + \alpha\zeta^2 \,d^2x\right]\,dt 
	\\
	&= \int\int \frac{D\rho}{2}|\bs{u}|^2 - p(D-1) + w(\p_t + \mathcal{L}_u)\zeta - \frac12(w^2 + \alpha\zeta^2) \,d^2x\,dt \,,
\end{aligned}
\end{equation}
where the areal density element, $Dd^2x$, and thermal buoyancy, $\rho$, are advected quantities. Taking variational derivatives of this with respect to $w$ and $\zeta$ gives, respectively, the following equations
\begin{align}
	(\p_t + \mathcal{L}_u)\zeta &= w
	\,,\\
	(\p_t + \mathcal{L}_u)w &= -\alpha\zeta
	\,.
\end{align}
The Euler-Poincar\'e equation is
\begin{equation}
\begin{aligned}
	(\p_t + \mathcal{L}_u)(D\rho\bs{u}\cdot d\bs{x} + wd\zeta) &= Dd\left( \frac{\rho}{2}|\bs{u}|^2 - p\right) - \frac{D}{2}|\bs{u}|^2 d\rho 
	\\
	&= D\rho\left( \frac{1}{2}|\bs{u}|^2\right) - Ddp \,.
\end{aligned}
\end{equation}
Putting together the equations for our harmonic oscillations, we have
\begin{equation}\label{eqn:HO_wave_momentum_advection}
	(\p_t+\mathcal{L}_u)(wd\zeta) = \frac12d\left( w^2 - \zeta^2 \right)\,,
\end{equation}
and hence our equations of motion are
\begin{align}
	\p_t\bs{u} + \bs{u}\cdot\nabla\bs{u} + \frac{1}{\rho}\Big(w\nabla w - \zeta\nabla\zeta\big) &= - \frac{1}{\rho}\nabla p \label{eqn:HO_EulerWaveEffect}
	\,,\\
	\p_t\rho + \bs{u}\cdot\nabla\rho &= 0
	\,,\\
	\nabla\cdot\bs{u} &= 0
	\,,\\
	\p_t\zeta + \bs{u}\cdot\nabla\zeta &= w
	\,,\\
	\p_tw + \bs{u}\cdot\nabla w &= -\alpha\zeta
	\,.
\end{align}
The effect of the current on the waves is that they now occur within a moving frame of reference, and the effect of the waves on the currents is given by the two new terms on the left hand side of the Euler momentum equation \eqref{eqn:HO_EulerWaveEffect}. Note that this interaction is a result of the addition of the term $\bs{u}\cdot(w\nabla\zeta)$ into the action.
\begin{remark}
	Within the Lagrangian we see a term of the form $w(\p_t +\mathcal{L}_u)\zeta$, in this instance this is equivalent to using a Lagrange multiplier to constrain the kinematic boundary condition by replacing this term with $\lambda(\p_t\zeta + \mathcal{L}_u\zeta - w)$.
\end{remark}
\begin{remark}
	The wave effect on current terms, arising due to equation \eqref{eqn:HO_wave_momentum_advection}, are absorbed into the pressure term. That is, we may write equation \eqref{eqn:HO_EulerWaveEffect} as
\begin{equation}
	\p_t\bs{u} + \bs{u}\cdot\nabla\bs{u} = - \frac{1}{\rho}\nabla\left(p +\frac{w^2}{2} - \frac{\zeta^2}{2} \right) \,,
\end{equation}
	and hence the potential vorticity form of the Euler equations is unchanged, because the harmonic oscillations simply contribute an additional term that redefines the definition of the pressure, but does not influence the velocity of the mean flow.
\end{remark}

As illustrated from this example, the coupling of harmonic oscillations to an incompressible Euler fluid, as appeared in \cite{Proceedings2022}, can be improved upon by a more physically relevant choice of wave Hamiltonian or fluid Lagrangian.

\end{appendices}

\end{document}